\documentclass[review,12pt]{elsarticle}

\usepackage{amssymb}
\usepackage{amsthm}
\usepackage{amsmath}
\usepackage{cases}
\usepackage{hyperref}

\usepackage{color,xcolor}
\usepackage{amsmath,amsfonts,amsthm,bm} 
\usepackage{subfigure} 
\usepackage{nomencl}
\usepackage{algorithm}
\usepackage{algorithmic}

\makenomenclature 

\usepackage{etoolbox}
\renewcommand{\nomgroup}[1]{%
	\item[\textbf{%
		\ifthenelse{\equal{#1}{D}}{Sets,Indices}{}
		\ifthenelse{\equal{#1}{P}}{Parameters}{}
		\ifthenelse{\equal{#1}{V}}{Variables}{}%
	}]%
}

\begin{document}

\begin{frontmatter}



\title{Plant-wide byproduct gas distribution under uncertainty in iron and steel industry via quantile forecasting and robust optimization}


%

\author[1]{Sheng-Long Jiang\corref{cor1}}%
\author[2]{Meihong Wang}
\author[3]{I. David L. Bogle}

\address[1] {College of Materials Science and Engineering, Chongqing University, China}
\address[2]{Department of Chemical and Biological Engineering, The University of Sheffield, U.K.}
\address[3] {Department of Chemical Engineering, University College London, U.K.}


\cortext[cor1]{Corresponding author: sh.l.jiang.jx@gmail.com}
\begin{abstract}
	Byproduct gas is one of the most important energy resources of the modern iron and steel industry because it  is crucial in supplying energy to manufacturing processes as well as converting energies, such as stream and electricity. The optimal distribution of byproduct gases in an iron and steel plant can significantly reduce energy costs and carbon emissions. However, the balance between supply and demand is easily threatened by the quantity- and quality-related uncertainties from manufacturing processes. Following the supply-storage-conversion-demand network, this study developed an optimal gas distribution model considering uncertain supply and proposed a two-stage robust optimization (TSRO) model including  ``here-and-now" decisions,  which minimize the start-stop cost of conversion units, and making ``wait-and-see" decisions, which minimize the operating costs of gasholders and demand penalties. To implement the TSRO model in practice, this study proposes a ``first quantify, then optimize" method: (1) quantify the uncertainty of surplus gas via a quantile regression-based multi-step time series model, and (2) find the optimal solution via a column-and-constraint generation algorithm. Finally, this study provides a case study of an industrial energy system to validate the proposed methodology.
\end{abstract}

%


\end{frontmatter}

\section{Introduction}
\label{S1:intro}

Given that the global population and living standards are improving, steel demand is expected to continuously grow in the coming decades, especially in developing countries, such as China, India, and Brazil \cite{Gahm2016}. However, the iron and steel industry is also an energy-intensive and high-pollution industry, contributing roughly $8\%$ of energy-related consumption and $6\%$ of global $\mathrm{CO}_2$ emission \cite{Fan2021}. This is a core challenge that the world must face if it is to achieve a low-carbon and sustainable manufacturing future. China is the largest steel producer worldwide with a crude steel production volume of more than one billion tons in 2021 and accounts for about 15\% of the total greenhouse gas emissions \cite{Ren2021}.

\begin{figure*}[htbp!]
	\centering
	\includegraphics[width=1.0\textwidth]{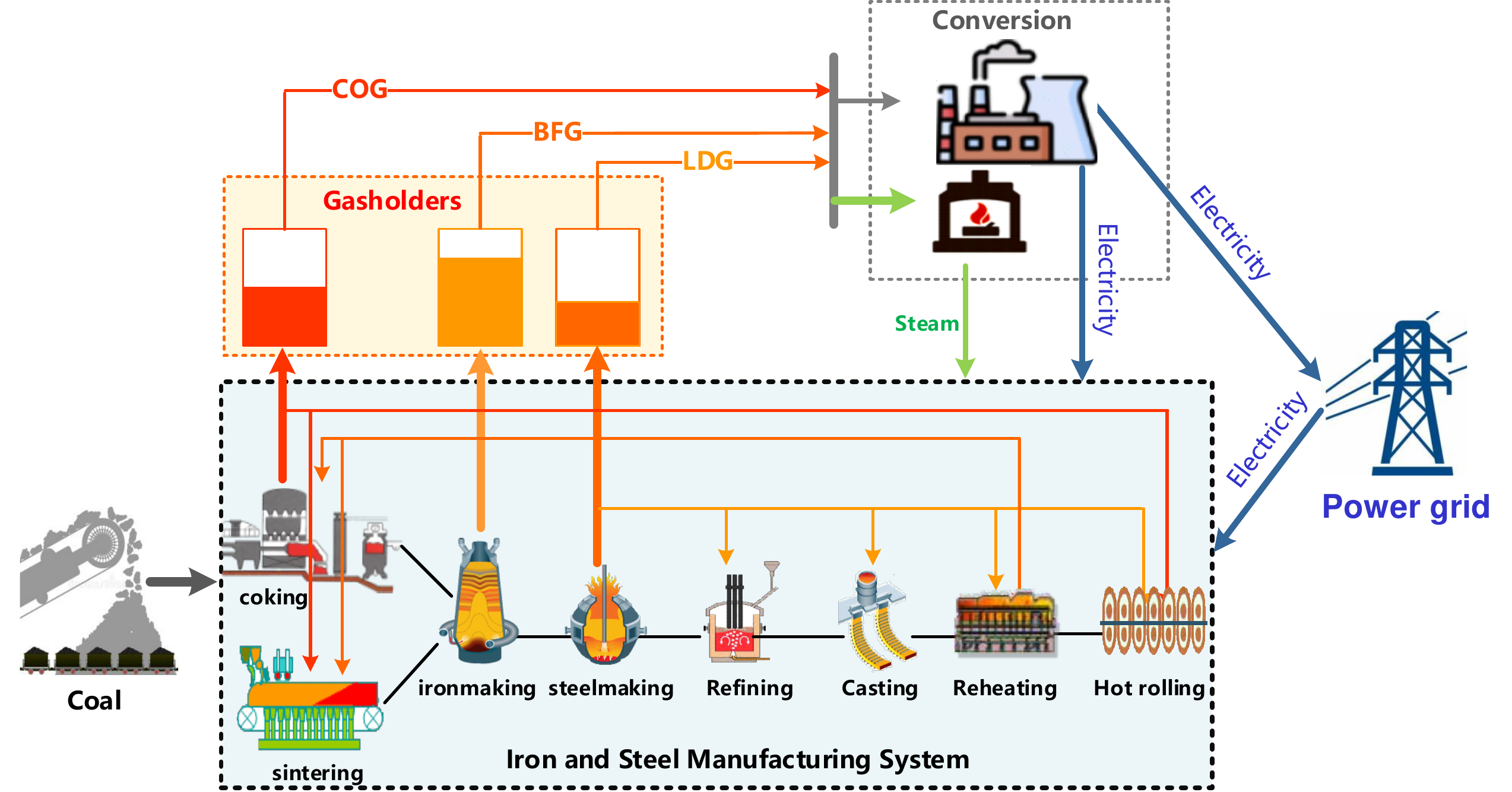}
	\caption{Byproduct gas system in a integrated iron and steel plant.}
	\label{fig:flowsheet}
\end{figure*}

In an integrated iron and steel plant, the iron ore to steel products primarily relies on coal-related resources and produces a large volume of byproduct gas, which is identified as the most important secondary energy source since its consumption is approximately $30\%$ in the plant-wide energy system \cite{Zhao2017a}. The byproduct gases come from coke ovens, blast furnaces, and basic oxygen furnaces and are known as coke oven gas (COG), blast furnace gas (BFG), and Linze-Donawitz gas (LDG) respectively. As shown in Figure \ref{fig:flowsheet}, the byproduct gases are first delivered to other production units in the manufacturing system. Further, the surplus gases are stored in the dedicated gasholders, or converted to other energies (i.e. stream and power) to meet the demands of some manufacturing units. Generating energy from the byproduct gases can reduce procurement costs in the marketplace. If the surplus gases cannot realize the conversion and exceed the capacity of gasholders, they have to be emitted into the air by flaring, resulting in gratuitous energy loss and environmental pollution. By summarizing these processes, the units are divided into three parts:   
\begin{enumerate} [1)]
	\item The manufacturing system, which includes production units, serves as both supplies and customers.
	\item The storage system, which includes gasholders, serves as buffers. 
	\item The conversion system, which includes boilers, CHP, and CDQ units, serves as adjusters. 
\end{enumerate}

Therefore, how to make optimal decisions on distributing byproduct gas to these sub-systems is a vital path to improving energy efficiency and reducing emissions. However, several factors contribute to the uncertainty of the energy system. These include the following:

\begin{enumerate}[1)]
	\item Quality-related uncertainty, which means the calories of byproduct gases from the manufacturing system vary, because the combustible components (H2, CH4, CO) in gases change randomly.
	\item Quantity-related uncertainty, which means the volume of byproduct gases cannot be precisely estimated, because the process units in the manufacturing system are not fully controllable.
\end{enumerate}

Given that these factors are quite intricate, traditional optimal distribution models with perfect assumptions for byproduct gases tend to fail in practice. To overcome these shortages, this study presents a robust optimization (RO) method that can provide several flexible strategies to absorb uncertainty and help the byproduct gas system operate with low risks.  The contributions of this study can be summarized as follows:

\begin{enumerate}[1)]
	\item Formulate a network flow-based mathematical model to achieve optimal distribution of byproduct gases in integrated iron and steel plants.   
	\item Develop a two-stage RO (TSRO) model to make "here-and-now" and "wait-and-see" decisions under uncertainty.
	\item Propose a solution methodology to the TSRO with the "first quantification, then optimization" idea.
\end{enumerate}

The remainder of this study is organized as follows. Section \ref{S2:review} reviews the optimal distribution models in integrated iron and steel plants and the RO technique under uncertainty. Section \ref{S3:model} proposes an optimization-based mathematical programming model for the byproduct gas distribution and formulate it as a TSRO model under uncertainty. Section \ref{S4:DDRO} reformulates the TSRO to a data-driven two-stage RO (DD-TSRO) model following the idea of “first quantify, then optimize”. Section \ref{S5:case} presents computational studies that verify the effectiveness of the proposed approach. Section \ref{S6:concl}, presents the conclusions of the study and makes recommendations for future research.

\section{Literature review}
\label{S2:review}
The optimization-based mathematical model provides a powerful tool for making distribution decisions on byproduct gas in iron and steel plants. As early as 1980, Markland \cite{Markland1980} first constructed an optimal distribution model of byproduct gas via linear programming (LP) to improve energy efficiency and reduce the cost of purchased fuel. In 1991, Akimoto et al. \cite{Akimoto1991} initially introduced mixed-integer LP (MILP)-based distribution model which involves the operational cost of gasholders. Kim et al. built MILP-based models to determine the optimal gas distribution to simultaneously optimize multiple conflict objectives \cite{Kim2003a}, and let it meet the varying energy demands and prices \cite{Kim2003b}. Afterward, attention moved towards building optimal distribution models of byproduct gas, with MILP becoming the most popular modelling tool. For example, Kong et al. \cite{Kong2010} simultaneously optimized the distribution of byproduct gases in the storage, conversion, and manufacturing system. Zhao et al. investigated how to make optimal gas distributions with different objective factors \cite{Zhao2015} (i.e. boiler and gasholder penalty) and under time-of-use electricity price \cite{Zhao2017b}. Zeng et al. \cite{Zeng2018} assumed generation rates of byproduct gases were time-varied, which provided flexibility on the supply side. Hu et al. \cite{Hu2022} provided optimal gas distribution in the cogeneration system considering a new start-stop behaviour of energy conversion equipment. These studies provide insight into the decision process of the byproduct gas distribution. However, they are thought to be Utopia because they do not consider uncertainty and rely on fully known parameters.

In recent years, some practitioners have recommended integrating uncertainty into byproduct gas scheduling practice. Zhao et al. \cite{Zhao2016} proposed a data-driven two-stage optimization model where, at the prediction stage, a Gaussian kernel-based regression model is used to predict the consumption of the outsourcing energy and at the second stage, a mathematical programming model sought the optimal solutions in a moving horizon way. Jin et al. \cite{Jin2018} established a causal reasoning model to predict the interval of each gasholder and developed a four-layer causal network to construct candidate solutions as well as select the best solution with an evaluation indicator. Pena et al.\cite{Pena2019} proposed different time series prediction models considering the difference between continuous and discrete processes, after which a MILP model is used to construct a dynamic regulation model in a moving horizon way. Focusing on the demand for carbon capture, utilization, and storage in iron and steel enterprises, Xi et al. \cite{Xi2021} used a gradient lifting regression tree as the proxy model of the energy conversion system and then used the PSO algorithm to find the optimal distribution decision. These studies all assumed the decision-maker has limited information but can accurately forecast critical parameters via machine learning techniques.

However, the uncertainty in the real-world byproduct gas system is highly intricate, and little information is known. RO is a novel technique that can work with little information about the underlying uncertainty, except for its upper and lower bounds, and seeks an optimal solution that covers the worst-case cost within a well-synthesized ``uncertainty set”. A comprehensive summary of developing and applying RO can be found in the tutorial \cite{Bertsimas2011} and the reviews presented by Gabrel et al. \cite{Gabrel2014} and Rahimian et al. \cite{Rahimian2019}. Little information in classic RO always makes the optimal solution over-conservative, with many practitioners nowadays trying to discover more information from historical data to develop a new style of RO, called data-driven or distributionally RO (DRO) \cite{Ning2019}. Delage et al. \cite{Delage2010} developed moment-based DRO model with mean and variance information. Bertsimas et al. \cite{Bertsimas2018} leveraged statistical hypothesis testing to find a tight value-at-risk bound and calibrate uncertainty sets. Shang et al. \cite{Shang2017} and Ning et al. \cite{Ning2018} derived data-driven polyhedron uncertainty sets using various machine learning techniques, such as support vector machine (SVM), principal component analysis (PCA), and kernel density estimation (KDE). DRO has also been applied to energy management systems (EMS) in industrial processes \cite{Qiu2022}. Zhao et al. \cite{Zhao2019} proposed to the KDE-based DRO model to address the operational optimization problem of industrial steam systems under uncertainty. Shen et al. \cite{Shen2020} applied an SVM-based DRO model in the energy system of an ethylene plant. 

Recent studies tend to focus on the uncertain energy distribution problem in integrated iron and steel plants and assumed that the key parameters (i.e. gasholder level) or the relationship between supply and demand under uncertainty were predictable and found optimal decisions with deterministic optimization models. However, a large variety of unforeseen events in realistic manufacturing and energy systems may affect the prediction accuracy and the feasibility of optimization models. To construct a more legitimate model, this study modelled the byproduct gas distribution problem with uncertain gas supply from the manufacturing system and calibrated the uncertainty set via machine learning.

\section{Optimal distribution under supply uncertainty}
\label{S3:model}

\subsection{Network of byproduct gas system}

\begin{figure*}[htbp!]
	\centering
	\includegraphics[width=0.9\textwidth]{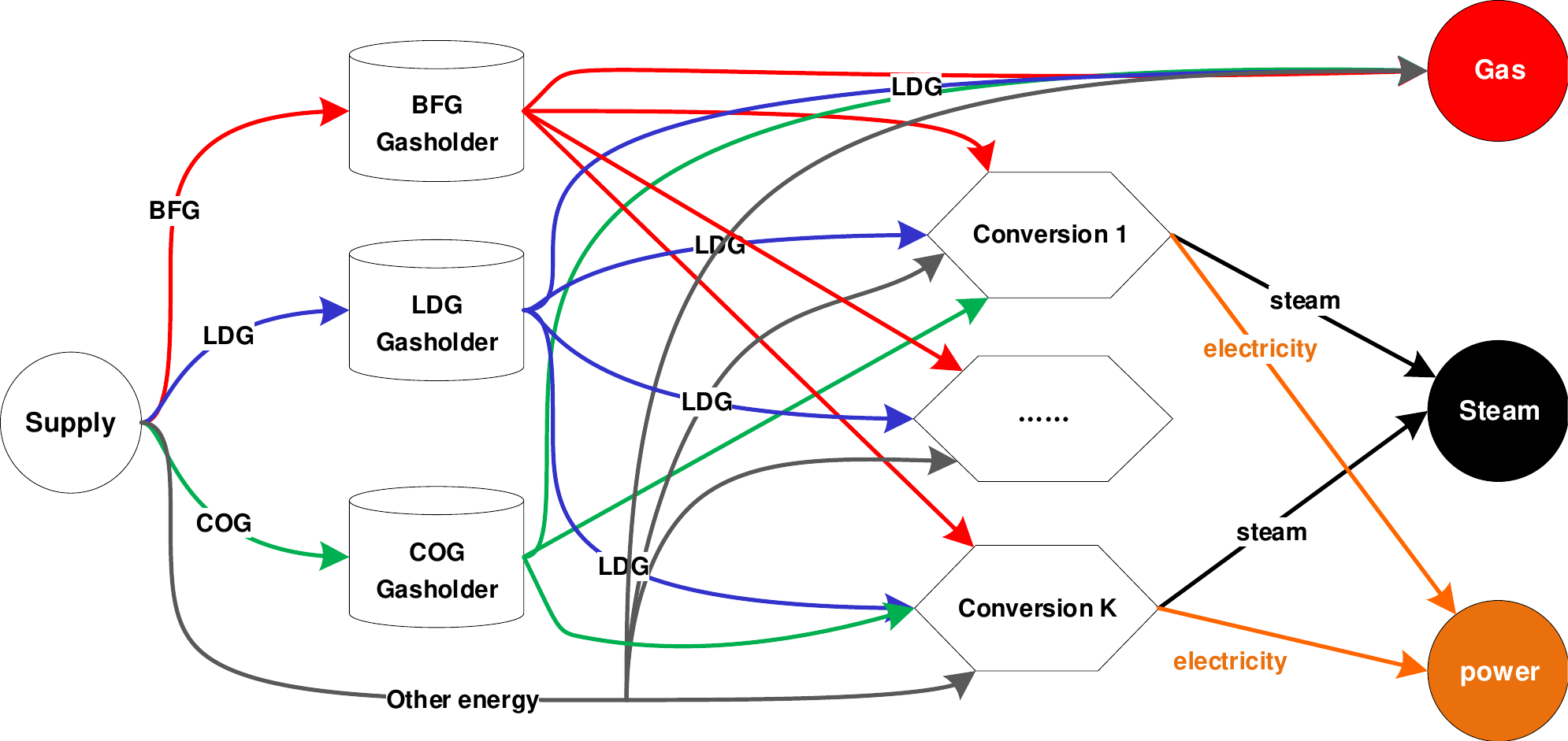}
	\caption{Network flow of byproduct gas system.}
	\label{fig:network}
\end{figure*}

This study considers the network of byproduct gas system ($\mathcal{K, A}$), where $\mathcal{K}$ denotes processing units and $\mathcal{A}$ denotes the flows between units. All processing units consist of the following four roles: 
\begin{enumerate}[1)]
	\item Supply ($\mathcal{K^+}$), which supplies surplus gases and other energies to the system.  
	\item  Storage ($\mathcal{K^S}$), which includes the gasholders of BFG, COG and LDG.
	\item  Conversion ($\mathcal{K^C}$), which converts the surplus gases to other energies, i.e electricity, and streams.
	\item Demand ($\mathcal{K^-}$), which denotes the energy demands of the iron and steel plant.
\end{enumerate}

Assumed that the EMS collects all required information, the mathematical model make the optimal decisions for byproduct gases distribution over a finite time horizon $\mathcal{T} = \{1,2,\dots,t,\dots, T\}$, where $t$ is the index of a period. 

\subsection{Deterministic optimization model}

\begin{enumerate}[1)]
	
\item Storage units: 
  
\quad\; In time period $t$, the byproduct gases flow into their dedicated gasholders ($k \in \mathcal{K}^{S}$) and the varying gasholder level satisfies mass balance \cite{Zeng2013}. This is the difference between period $t-1$ and $t$ is equal to the surplus volume of byproduct gas minus the total volume of the output flows and can be described by the following linear equation:
\begin{equation}
	\label{eq:Vol}
	u_{k,t} - u_{k,t-1} = \sum_{a \in \mathcal{A}^{+}(k)} z_{a,t}-\sum_{a \in \mathcal{A}^{-}(k)}f_{a,t}, \quad k \in \mathcal{K}^{S} , t \in \mathcal{T}
\end{equation}
where $u_{k,0}$ is the initial gasholder level. 

\quad\; Since each gasholder must operate within a safe operational region, its level has the following lower and upper bounds:
\begin{equation}
	\label{eq:Vbnd}
	    \underline{U}_{k}	\leq  u_{k,t}, \leq \overline{U}_{k}, \quad k \in \mathcal{K}^{S}, t \in \mathcal{T}
\end{equation}

\quad\; Due to operational restraints, the gasholder level change $|u_{k,t} - u_{k,t-1}|$ of gasholder $k$ between periods $t-1$ and $t$ must not excess $\Delta_{k}$.
\begin{equation}
	\label{eq:Vchg}		
	u_{k,t} - u_{k,t-1} \leq \Delta_{k}, u_{k,t-1} - u_{k,t} \leq \Delta_{k},  \quad k \in \mathcal{K}^{S} , t \in \mathcal{T}	
\end{equation}

\quad\; The deviation $v_{k,t}$ of gasholder $k$ from the middle level $u_{k,\rm{mid}}$ to the current level $u_{k,t}$ indicates the risk of under-stock or over-stock. Thus, the deviation ($v_{k,t}$) can be defined with the following inequations:
\begin{equation}
	\label{eq:Vdev}
	u_{k,t} - U_{k,\rm{mid}} \leq v_{k,t},   U_{k,\rm{mid}} - u_{k,t} \leq v_{k,t}, \quad k \in \mathcal{K}^{S} , t \in \mathcal{T}
\end{equation}

\item Conversion units

\quad\; Via a conversion unit (such as boiler and CHP), input flows can be converted into other output flows, which must satisfy the following input-output balance of energy conversion. 
\begin{equation}
	\label{eq:conv}
	\rho_{k} \sum_{a \in \mathcal{A}^{+}(k)} f_{a,t} \times \omega_{e:a} = \sum_{a \in \mathcal{A}^{-}(k)} f_{a,t} \times \omega_{e:a} , \quad k \in \mathcal{K}^C , t \in \mathcal{T}
\end{equation}

\quad\; When a conversion unit is turned on ($O_{k,t}=1$), its input and output energy flows are restricted by the capacity of conversion unit $k$.
\begin{equation}
	\label{eq:limt1}
	\underline{F}_{\ k,a}^{+} O_{k,t} \leq   f_{a,t}  \leq \overline{F}_{k,a}^{+} O_{k,t} , \quad k \in \mathcal{K}^C , a \in \mathcal{A}^{+}(k), t \in \mathcal{T}
\end{equation}
\begin{equation}
	\label{eq:limt2}
	\underline{F}_{\ k,a}^{-} O_{k,t} \leq   f_{a,t}  \leq \overline{F}_{k,a}^{-} O_{k,t} , \quad k \in \mathcal{K}^C , a \in \mathcal{A}^{-}(k), t \in \mathcal{T}
\end{equation}

\quad\; To ensure an equipment unit works at its normal conditions, the mixed calorific value of the input flows must be greater than the minimum value $\eta_k^+$.
\begin{equation}
	\label{eq:in-limit}
	\sum_{a \in \mathcal{A}^{+}(k)} f_{a,t} \times \omega_{e:a} \geq \eta_k^+ \sum_{a \in \mathcal{A}^{+}(k)} f_{a,t} \times+ bigM(O_{k,t}-1), \quad k \in \mathcal{K}^C , t \in \mathcal{T}
\end{equation}
	
\quad\;  Given that the operation of large-capacity equipment with low output is thought to be uneconomic \cite{Hu2022}, the ratio of output flow to its maximum limits must not be allowed to be lower than the defined threshold value $\eta_k^-$. Otherwise, the conversion unit should be closed. 
\begin{equation}
	\label{eq:out-limit1}
	\sum_{a \in \mathcal{A}^{-}(k)} f_{a,t}  \geq \eta_k^- \sum_{a \in \mathcal{A}^{-}(k)} \overline{F}^{-}_{k,a} + bigM(O_{k,t}-1), \quad k \in \mathcal{K}^C , t \in \mathcal{T} \\
\end{equation}

\begin{equation}
	\label{eq:out-limit2}	
	\sum_{a \in \mathcal{A}^{-}(k)} f_{a,t}  \leq bigM \times O_{k,t},  \quad k \in \mathcal{K}^C , t \in \mathcal{T} \\	
\end{equation}

\quad\; If a conversion unit is being started or stopped ($S_{i,k}$), its on/off status ($O_{k,t}$) is also being changed between periods $t-1$ and $t$. In a real-world situation, conversion units should stably operate and avoid repeated start-stop changes. Thus, the binary variables $S_{k,t}$ and $O_{k,t}$ are related as follows:
\begin{equation}
	\label{eq:start-stop}
	 O_{k,t} -  O_{k,t-1} \leq S_{k,t} , O_{k,t-1} - O_{k,t} \leq S_{k,t} \quad k \in \mathcal{K}^S , t \in \mathcal{T}
\end{equation}

\item Demand side

\quad\; The demand side has two types of energy: the emitted gases ($\mathcal{K}^{-}_{1}$), including BFG, LDG, and COG; and the produced energies ($\mathcal{K}^{-}_{2}$), including electricity and steam. In principle, the demands of byproduct gases (BDG, COG, LDG) are zero because a surplus will cause environmental pollution and shortage will leads the safety risk of gasholders. The surplus of other energies means the profits on sale, whereas the shortage will result in a purchase cost. The total generating energy (input of demand $k$) at period $t$ must meet its demand $d_{k,t}$.

\begin{equation}
	\label{eq:dmd}
	\left\{
	\begin{split}		
		\sum_{a \in \mathcal{A}^{+}(k)} f_{a,t} - w_{k,t} \geq d_{k,t} ,& \quad k \in \mathcal{K}^{-,1}, t \in \mathcal{T} \\		
		\sum_{a \in \mathcal{A}^{+}(k)} f_{a,t} + w_{k,t} \geq d_{k,t} ,& \quad k \in \mathcal{K}^{-,2}, t \in \mathcal{T}.		
	\end{split}	
	\right.
\end{equation}
where $w_{k,t}$ indicates the surplus and shortage of demands in $\mathcal{K}^{-}_{1}$ and $\mathcal{K}^{-}_{2}$. 

\item Objectives

\quad\; Given the assumptions and constraints above, the optimal objective of the studied byproduct gas system within the time horizon ($\mathcal{T}$) can be represented as follows:   
\begin{equation}
	\label{eq:objfunc}
	\min f =   \gamma_1 \sum_{t \in \mathcal{T}} \sum_{k \in \mathcal{K}^C } S_{k,t}		
	+ \gamma_2\sum_{t \in \mathcal{T}} \sum_{k \in \mathcal{K}^C }  v_{k,t} 
	+  \sum_{t \in \mathcal{T}} \sum_{k \in \mathcal{K}^{-}} \gamma_{3,k}  w_{k,t} 
\end{equation}
where the first term represents the start-stop costs of conversion units, the second term denotes the operating costs of gasholder deviation from the middle position, and the last term indicates the surplus or shortage cost of energy demands. It should be noted that in \eqref{eq:objfunc}, $\gamma_1$ and $\gamma_2$ are unit-unrelated constant coefficients, and $\gamma_{3,k}$ is a unit-related coefficient. Elsewhere, $\gamma_2$ is represented by a set of piece-wised functions \cite{Zhao2015}. 

\end{enumerate}

\subsection{TSRO model}

It should be noted that, in the gasholder level balance equation \eqref{eq:Vol}, the volume of gas supply ($z_{a,t}$) is equal to the generation of manufacturing units minus their consumption \cite{Zhao2017a}. Given that manufacturing units are not fully controlled, the gas volumes are random and equation \eqref{eq:Vol} needs to hold under the random. In this study, it was assumed that $z_{a,t}$ is an uncertain variable that varies between $[z_{a,t}^{\circ}-\hat{z}^{-}_{a,t}, z_{a,t}^{\circ} + \hat{z}^{+}_{a,t}]$, where $z_{a,t}^{\circ}$ is the forecasting value (i.e. nominal value) and $\hat{z}^{-}_{a,t}, \hat{z}^{+}_{a,t}$ respectively denote its maximum negative and positive forecasting deviations. Thus, the random supplies over time can be represented as the following box-uncertainty set.

\begin{equation}
	\label{eq:uset}
	\mathcal{Z}_{\rm{box}} := \left\lbrace \bm{z}: z_{a,t} = z_{a,t}^{\circ} + \xi_{a,t}^{+}\hat{z}_{a,t}^+ - \xi_{a,t}^{-}\hat{z}_{a,t}^{-}, 0 \leq \xi_{a,t}^{-},\xi_{a,t}^{+} \leq 1, \forall a, t \right\rbrace 
\end{equation} 

Since the box-uncertainty set always suffers from over-conservativeness, Bertsimas and Thiele \cite{Bertsimas2006} introduced an integer parameter  $\Gamma_t(0 \leq \Gamma_t \leq T)$ (budget of uncertainty) to restrict the maximum cumulative deviation with the summation of absolute value constraints.

\begin{equation}
	\label{eq:bset}
	\begin{aligned}
	\mathcal{Z}_{\rm{bud}} := & \left\lbrace \bm{z}: z_{a,t} = z_{a,t}^{\circ} + \xi_{a,t}^{+}\hat{z}_{a,t}^+ - \xi_{a,t}^{-}\hat{z}_{a,t}^{-}, \right. \\ 
	 & \left.  0 \leq \xi_{a,t}^{-},\xi_{a,t}^{+} \leq 1, \sum_{t=1}^{T} \left(\xi_{a,t}^{+} + \xi_{a,t}^{-} \right) \leq \Gamma_a, \forall a,t \right\rbrace 
	\end{aligned}
\end{equation} 

With the defined uncertainty sets above, a TSRO model can be formulated to make the decisions of which units are started or stopped in the first-stage minimization, whilst the actual operating variables can be determined in the second-stage minimization after the potential variations have been realized via the maximization over uncertainty set $\mathcal{Z}$. 

\begin{equation}
	\label{eq:TSRO}
	\begin{split}
	\min_{\substack{O, S \in \{0,1\}}} \left\lbrace  \gamma_1 \sum_{t \in \mathcal{T}} \sum_{k \in \mathcal{K}^C }  S_{k,t} 			
	 + \max_{z \in \mathcal{Z} } \min_{{\substack{u,v,w,f \in \mathbb{R}}}}   \sum_{t \in \mathcal{T}}  \left(\sum_{k \in \mathcal{K}^S }  \gamma_2 v_{k,t} + \sum_{k \in \mathcal{K}^- } \gamma_{3,k}  \times w_{k,t}\right) \right\rbrace \\
	\end{split}
\end{equation}
subject to constraints \eqref{eq:Vol}-\eqref{eq:dmd}

\section{DDRO formulation and solving}
\label{S4:DDRO}

\begin{figure*}[htbp!]
	\centering
	\includegraphics[width=0.7\textwidth]{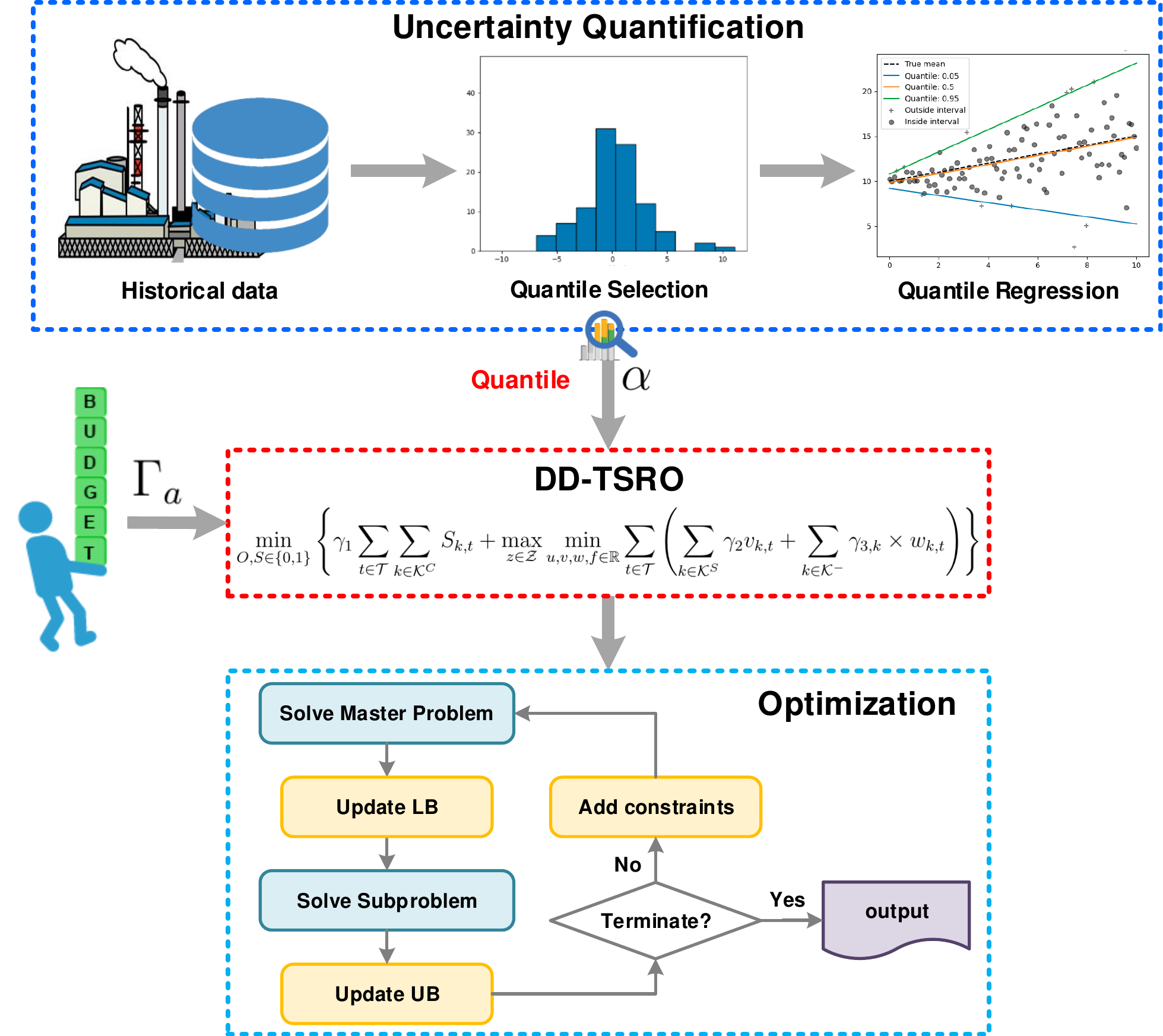}
	\caption{Flowchart of the proposed DD-TSRO framework.}
	\label{fig:solution}
\end{figure*}

The formulations presented in the previous section offer a complete framework to apply TSRO to byproduct gas distribution models, however, the mathematical model needs to link to the practice. In this section, a reformulation the TSRO to a DD-TSRO model following the idea of “first quantify, then optimize” is given, with detailed steps illustrated in Figure \ref{fig:solution}

\subsection{Uncertainty quantification via machine learning}

Since the RO technique is interesting in forecasting intervals instead of points, the conditional quantile regression (CQR) proposed by Koenker et al. \cite{Koenker1978} can be applied to define a time series model. The CQR aims to quantify the uncertain gas supply at a period $\tau$ (${z}_{a,\tau}$ where subscript $a$ is omitted in this section) given data of past $p$ periods. 
\begin{equation}
	\label{eq:ts}
	 \tilde{z}_{\tau} = \psi(z_{\tau-1},\dots,z_{\tau-p};\theta)+\epsilon_{\tau}
\end{equation}
where $\tilde{z}_{\tau}$ represents the prediction, $\psi(\cdot)$ represents the latent function of unknown form, $\theta$ represents its hyper-parameter, and $\epsilon_{\tau}$ represents the random noise, but is not necessarily normally distributed. Instead of forecasting the conditional mean $\mu(\tilde{z}_{\tau}|z_{\tau-1},\dots,z_{\tau-p})$, the CQR technique can be applied to estimate a given quantile of $\tilde{z}_{\tau}$, such as the median. Set $L_{\alpha}$ as the $\alpha$th quantile of the cumulative density function of $\tilde{z}_{\tau}$, i.e. ${\rm Pr}(\tilde{z}_{\tau} \leq L_{\alpha})=\alpha$ and the $\alpha$th conditional quantile function of ${z}_{\tau}$ may be rewritten as:
\begin{equation}
	\label{eq:QF}
	q_{\alpha}(\tilde{z}_{\tau}|z_{\tau-1},\dots,z_{\tau-p}) := \inf\{{z}_t \in \mathbb{R}: {\rm Pr}(\tilde{z}_{\tau}|z_{\tau-1},\dots,z_{\tau -p})\geq \alpha\}
\end{equation}

Given the training set $\mathcal{D} = \{ (z_{\tau-1}^{(i)},\dots,z_{\tau-p}^{(i)}; z_{\tau}^{(i)})\}_{i=1}^{N}$, classical regression analysis estimates the conditional mean by minimizing the sum of squared residuals on training sets.  Analogously, CQR estimates a conditional quantile function $q_{\alpha}$ of ${z}_{\tau}$ by minimizing the following loss function:
\begin{equation}
	\label{eq:loss}
	\mathcal{L} = \min_{\theta} {\frac{1}{N} \sum_{i=1}^{N}\left[ \alpha/2 \max(z_{\tau}^{(i)} -\tilde{z}_{\tau}^{(i)},0) + (1-\alpha/2) \max(z_{\tau}^{(i)} -\tilde{z}_{\tau}^{(i)},0)\right] }
\end{equation}

When $\alpha < 0.5$, smaller prediction garners more attention in loss; when  $\alpha > 0.5$, larger prediction garners more attention in loss; with $\alpha = 0.5$, it reduced to least absolute deviation regression (median regression). Therefore, given the training set  $\mathcal{D}$ and quantile $\alpha$, the uncertainty of $z_t$can now be quantified:
\begin{equation}
	\label{eq:intval}
	\begin{cases}
		{z}_{\tau}^{\circ}= q_{0.5} (\tilde{z}_{\tau}|{z}_{\tau-1},\dots,{z}_{\tau-p})\\
		\hat{z}_{\tau}^{-}= {z}_{\tau}^{\circ} - q_{\alpha} (\tilde{z}_{\tau}|{z}_{\tau-1},\dots,{z}_{\tau-p}) \\
		\hat{z}_{\tau}^{+}= q_{1-\alpha} (\tilde{z}_{\tau}|{z}_{\tau-1},\dots,{z}_{\tau-p}) -{z}_{\tau}^{\circ} 
	\end{cases}
\\
\end{equation}

The gradient boosting decision tree (GBDT) is applied to find the best latent function $\psi(\cdot)$ to cover the training data. The GBDT is a powerful machine learning technique with multiple regression trees (e.g. CART). The prediction model can be rewritten as follows:
\begin{equation}
\psi_M({z}_{\tau-1},\dots,{z}_{\tau-p})= \sum_{m=1}^{M} tree ({z}_{\tau-1},\dots,{z}_{\tau-p};\Theta_m) 
\end{equation}
where $tree(\cdot ; \cdot)$ represents a regression tree, $\Theta_m$ denotes parameters of tree $m$, and $M$ represents the number of regression trees. Then,  the GBDT uses a gradient boosting and a forward step-wise algorithm to find the optimal hyper-parameters. 
\begin{equation}
	\Theta_m^*= \arg \min_{\Theta_m} \mathcal{L} (\tilde{z}_{\tau}, \psi_{m-1}+ tree ({z}_{\tau-1},\dots,{z}_{\tau-p};\Theta_m)
\end{equation}

The details refer to the seminal work of \cite{Friedman2001}. This study also proposes a framework for implementing T-step prediction, as shown in the following:  
\begin{equation} \tilde{z}_{\tau+t} = \psi_M^{t}({z}_{\tau-1},\dots,{z}_{\tau-p},\alpha), t \in {1,....,T}  \end{equation}

\subsection{Reformulation for TSRO}

For formulation simplicity, the TSRO model is presented in the following matrix form:
\begin{subequations}
	\label{eq:STSRO}
	\begin{align}
	\min &\quad \mathbf{c}^\top \mathbf{x} + \max_{\mathbf{z} \in \mathcal{Z} }\min_{y \in \Omega(\mathbf{x,z})} \mathbf{d}^\top\mathbf{y} \label{eq:STSRO-1} \\
	s.t. &\quad \mathbf{A} \mathbf{x} \leq \mathbf{b}, \mathbf{x} \in \{0,1\} \label{eq:STSRO-2} \\
	 	 &\quad \Omega(\mathbf{x,z}) =
	 	 \begin{cases} 
	 	 	\mathbf{G} \mathbf{y} \geq \mathbf{h}  \\ 
	 	 	\mathbf{Q} \mathbf{y} \geq \mathbf{r} - \mathbf{P} \mathbf{x}   \\
	 	 	\mathbf{W} \mathbf{y} =  \mathbf{z} + \mathbf{s} \label{eq:STSRO-3}
 	  \end{cases} 	
   \end{align}
\end{subequations}

In formulation \eqref{eq:STSRO}, vector $\mathbf{x}$ represents the first-stage decisions ($O,S$), and vector $\mathbf{y}$ represents the second-stage decisions ($u,v,w,f$). The objective function \eqref{eq:STSRO-1} is divided into two parts: one depends on the binary variables  $\mathbf{x}$, and the other depends on the continuous variables  $\mathbf{y}$. Equation \eqref{eq:STSRO-2} includes all constraints involving only binary variables \eqref{eq:limt1},\eqref{eq:limt2}, \eqref{eq:out-limit2} and \eqref{eq:start-stop}.In the domain defined by equation \eqref{eq:STSRO-3}, the first term collects constraints  \eqref{eq:Vbnd}-\eqref{eq:conv} only involving continuous variables, the second term accounts for \eqref{eq:in-limit} and  \eqref{eq:out-limit1} that involving mixed variables, and the last term represents the constraints \eqref{eq:Vol} that involves the uncertain gas supply. 

With the simplified formulation, the dual of the second-stage optimization problem is first found to be LP and holds a strong duality. Therefore, its dual problem can be rewritten in the following form:
\begin{equation}
	\label{eq:sp1}
	\begin{split}
		 \max_{\mathbf{z},\bm{\lambda,\sigma,\phi}} &\quad \bm{\lambda}^\top \mathbf{h}+ \bm{\sigma}^\top (\mathbf{r-Px}) +  \bm{\phi}^\top (\mathbf{z+s})  \\
		s.t.  & \quad \bm{\lambda}^\top \mathbf{G}+ \bm{\sigma}^\top \mathbf{Q} +  \bm{\phi}^\top \mathbf{W}  = \mathbf{d}^\top\\
		& \quad \bm{\lambda} \geq 0 ,\bm{\sigma} \geq 0, \mathbf{z} \in \mathcal{Z}	
	\end{split}
\end{equation}
where $\bm{\lambda},\bm{\sigma},\bm{\phi}$ are the Lagrangian multipliers of the formula \eqref{eq:STSRO-3} and $\bm{\phi}$ is unbounded. Note that $\bm{\phi}^\top \mathbf{z}$ are bilinear in the objective function, and therefore needs to be linearized.  First, due to the variable $z$ being independent of other variables in the equation \eqref{eq:sp1}, the optimal solution must inferred to be one extreme point of $\mathcal{Z}$. Next, the uncertain set is injected into the term ( $\bm{\phi}^\top \mathbf{z}$) to obtain the following equation:
$$ \bm{\phi}^\top \mathbf{z}  = \sum_{t=1}^{T}\left( \phi_tz_t^\circ + \phi_t\xi_t^+\hat{z}_t^+ - \phi_t\xi_t^{-}\hat{z}_t^{-}  \right)  $$

Then, \eqref{eq:sp1} is transformed to the following MILP model using the Big-M method \cite{Guo2016}, which can be easily solved by commercial solvers. 
\begin{equation}
	\label{eq:sp2}
	\begin{split}
		\max_{\bm{\pi},\bm{\lambda,\mu,\phi}} &\quad \bm{\lambda}^\top \mathbf{h}+ \bm{\sigma}^\top (\mathbf{r-Px}) + \bm{\phi}^\top \mathbf{s} + \sum_{t=1}^{T}\left( \phi_t{z_t^\circ} + \pi_t^+\hat{z}_t^+ + \pi_t^{-}\hat{z}_t^{-}  \right)  \\
		s.t.  & \quad \bm{\lambda}^\top \mathbf{G}+ \bm{\sigma}^\top \mathbf{Q} +  \bm{\phi}^\top \mathbf{W}  = \mathbf{d}^\top\\
		& \quad \pi_{t}^{+} \leq {bigM} \xi_{t}^{+}, \pi_{t}^{+} \leq \phi_t + {bigM} \left( 1-\xi_{t}^{+}\right), \forall t \\
		& \quad \pi_{t}^{-} \leq {bigM} \xi_{t}^{-}, \pi_{t}^{-} \leq {bigM} \left( 1-\xi_{t}^{-}\right) - \phi_t, \forall t \\		
		& \quad \sum_{t=1}^{T} \left( \xi_{t}^+ + \xi_{t}^- \right) \leq \Gamma,  \xi_{t}^+, \xi_{t}^- \in \{0,1\}	\\
		& \quad \bm{\lambda} \geq 0 ,\bm{\sigma} \geq 0 \\
	\end{split}
\end{equation}
where $\pi_{t}^{+}$ and $\pi_{t}^{-}$ are auxiliary variables. If $\xi_{t}^{+}$ and $\xi_{t}^{-}$ are equal to 1, $\pi_{t}^{+}$ and $\pi_{t}^{-}$ will be limited to $\phi_t$ and $-\phi_t$. If $\xi_{t}^{+}$ and $\xi_{t}^{-}$ are equal to 0, $\pi_{t}^{+}$ and $\pi_{t}^{-}$ will be limited 0. 

\subsection{Column-and-constraint generation algorithm}
\begin{subequations}
	\label{eq:mp}
	\begin{align}
		\max_{\mathbf{x},\beta, \mathbf{y}_i|i<I} &\quad \mathbf{c}^\top \mathbf{x} + \beta \label{eq:CCG-1} \\
		s.t. &\quad \mathbf{A} \mathbf{x} \leq \mathbf{b}, \mathbf{x} \in \{0,1\} \label{eq:CCG-2} \\ 	
		& \quad \beta \geq \mathbf{d}^\top\mathbf{y}_i, i= 1,\dots,I \label{eq:CCG-3} \\
		& \quad \mathbf{G} \mathbf{y}_i \geq \mathbf{h} , i= 1,\dots,I \label{eq:CCG-4} \\  
		&\quad \mathbf{Q} \mathbf{y}_i \geq \mathbf{r} - \mathbf{P} \mathbf{x}  , i= 1,\dots,I \label{eq:CCG-5}\\ 
		&\quad \mathbf{W} \mathbf{y}_i =  \mathbf{z}_i + \mathbf{s} , i= 1,\dots,I, 	
	\end{align}
\end{subequations}

It is difficult to directly find the optimal solution of \eqref{eq:STSRO} within a short time via state-of-the-art optimization solvers (e.g. Cplex, Gurobi, and SCIP). Since the uncertainty set $\mathcal{Z}$ is a polyhedron, the worst scenarios of surplus gas are only located at an extreme point of $\mathcal{Z}$; therefore, the number of possible worst scenarios is finite \cite{Zeng2013}. According to this characteristic, this study applied the column-and-constraint generation (C\&CG) algorithm proposed by Zeng and Zhao \cite{Zeng2013} to seek the optimal solution for DD-TSRO. Initially, the TSRO defined in \eqref{eq:STSRO} is relaxed and reformulated as the master problem \eqref{eq:mp}. Then, the possible worst-case scenarios are identified by solving the subproblem \eqref{eq:sp2} and added into \eqref{eq:mp} as a cutting plane. Finally, the optimal solution can be found via iteratively improving the gap between the lower and upper bound. The solving procedure is stated in Algorithm~\ref{alg:CCG}.

\begin{algorithm}	
	\renewcommand{\algorithmicrequire}{\textbf{Input:}}
	\renewcommand{\algorithmicensure}{\textbf{Output:}}
	\caption{Column-and-Constraint Generation (C\&CG)}
	\label{alg:CCG}
	\begin{algorithmic}[1]
		\STATE Initialization: $ \rm{LB} = -\infty$, $ \rm{UB} = \infty$, $i = 1$
		\STATE Construct a feasible decision at the first stage $\mathbf{x}_1$   (via the deterministic model \eqref{eq:Vol}-\eqref{eq:objfunc})
		\REPEAT
		\STATE  Solve subproblem \eqref{eq:sp2} with fixed $\mathbf{x}_i$, obtain the optimal solution $\mathbf{\xi}'_i$ and the optimal objective value ${\beta}'$.  Then, a worst-case scenario $\mathbf{z}_i$ is construct via the budget-based set \eqref{eq:bset}.
		\STATE Update $\rm{UB}= \min \{\rm{UB}, \mathbf{c}^{\top}\mathbf{x}_i + {\beta}'\}$ 
		\STATE Solve the master problem defined in \eqref{eq:mp} by adding the new scenario $\mathbf{z}_i$  and new variables $\mathbf{y}_i$. Let $(\mathbf{x}^*, \beta^*)$ be the optimal solution. 
		\STATE  Update $ \rm{LB} = \mathbf{c}^{\top}\mathbf{x}^* + \beta^*$.
		\STATE  $i\leftarrow i + 1,  \mathbf{x}_i \leftarrow \mathbf{x}^*$.
		\UNTIL $\rm{UB}-\rm{LB} \leq \delta$ or $i>I$ ($\delta$ is a user-defined tolerance)  
		\ENSURE  The optimal solution $\mathbf{x}^*$, LB and UB.
	\end{algorithmic}  
\end{algorithm}

\section{Case study}

This section presents a case study on the energy system of an integrated iron and steel plant in China and verify the proposed approach. This steel plant has three gasholders, two 35t boilers, two 130 t boilers, two sets of CHP, and two sets of CDQ. The operational parameters and the network configuration of the studied energy system are presented in the attached supplementary materials (ref. Tables S1-S6)

For uncertainty quantification, the open-source scikit-learn package, version 1.1.1 \cite{Pedregosa2011} was applied, and the DD-TSRO model was programmed by Pyomo \cite{Hart2011}(a Python-coded, open-source modeling language) and solved by  “Gurobi 9.0” (\url{https://www.gurobi.com/products/gurobi-optimizer/}, an academic version with default settings).  All the computational studies were executed on a PC with an Intel Core i7 processor (3.60 GHz), 16 GB RAM, and a Windows 10 operating system.

\label{S5:case}

\subsection{Uncertainty quantification}

\begin{figure*}[htbp!]
	\centering
	\includegraphics[width=0.6\textwidth]{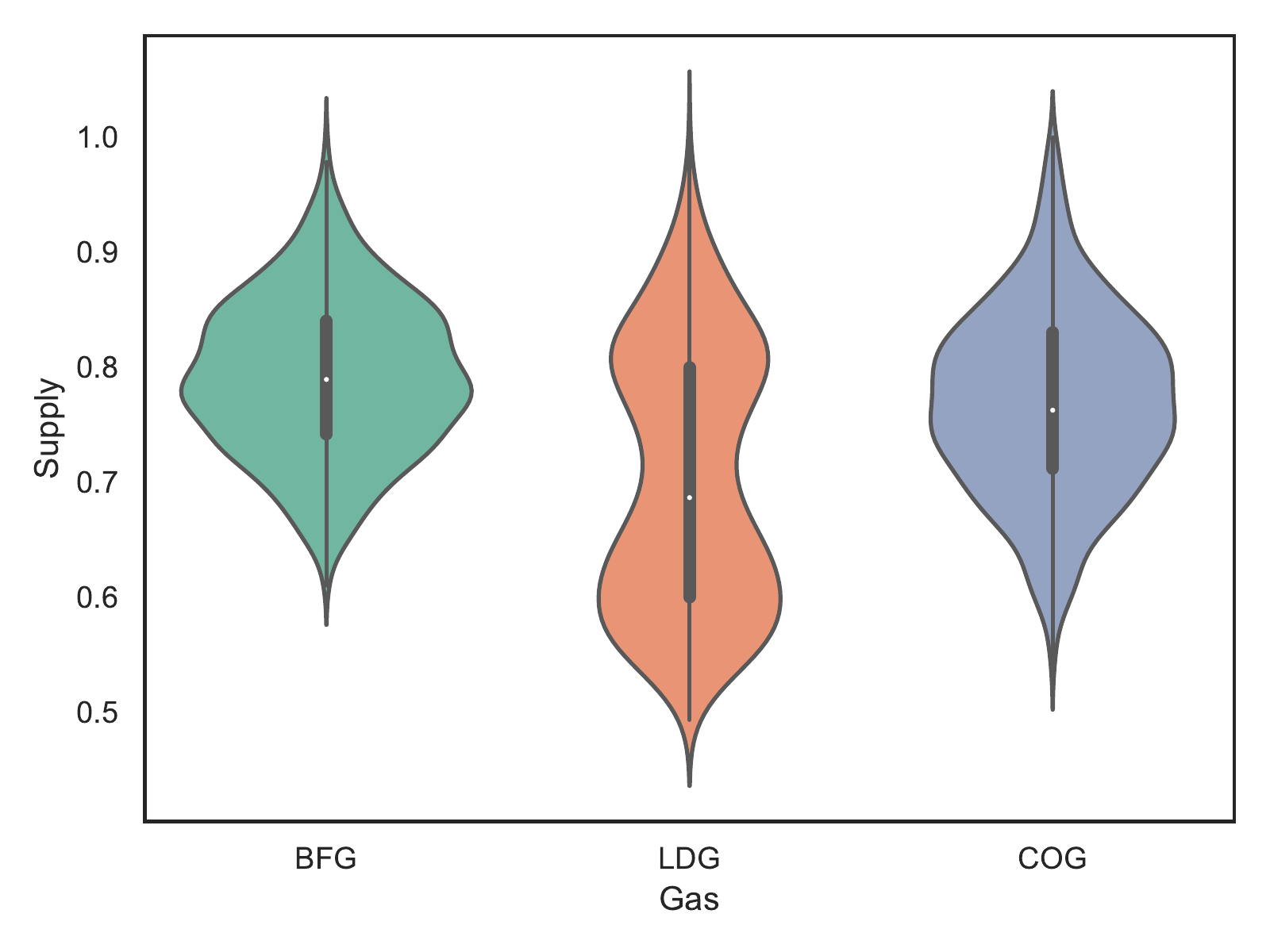}
	\caption{The violin plot of gas type and supply on dataset.}
	\label{fig:violin}
\end{figure*}

There are 1000 historical data points (ref. the attached data file is store at the webpage  \url{https://github.com/janason/Energy/tree/main/gas}) collected from the iron and steel plant. The GBDT technique was implemented to quantify the uncertain parameters and coded with scikit-learn's class via \textit{GradientBoostingRegressor}. More specifically, the dataset was separated into training (with 900 samples) and testing (with 100 samples) sets. The dataset was statistically analysed, with the violin plot of their normalizations shown in Figure \ref{fig:violin}. 

In the proposed T-step time-series model, $T = 8, p = 20$. When forecasting the nominal values, $\alpha = 0.5$ and when forecasting the deviations, the quantile $\alpha$ was selected from $\{0.01, 0.05, 0.1\}$. To observe the prediction effects, the performance of the quantification approach was evaluated in two aspects: (1) point forecasting metric, called mean absolute percentage error (MAPE), and (2) interval forecasting metric,  called prediction interval coverage probability (PICP):

$$ \mathrm{MAPE} = \frac{1}{N} \sum_{\tau=1}^{N} \left| \frac{z^{\circ}_{\tau}- \tilde{z}_{\tau}}{\tilde{z}_{\tau}} \right| \times 100\% $$

$$ \mathrm{PICP}= \frac{1}{N} \sum_{\tau=1}^{N} Cov_{\tau}^{(\alpha)} $$
where $$  Cov_{\tau}^{(\alpha)} = \begin{cases} 1, & \text{if} \quad z_{\tau} \in (z^{\circ}_{\tau}-\hat{z}^{-}_{\tau},z^{\circ}_{\tau}+\hat{z}^{+}_{\tau}) , \\ 0, & \text{otherwise} \end{cases}  $$

\begin{figure*}[htbp!]
	\centering
	\includegraphics[width=0.4\textwidth]{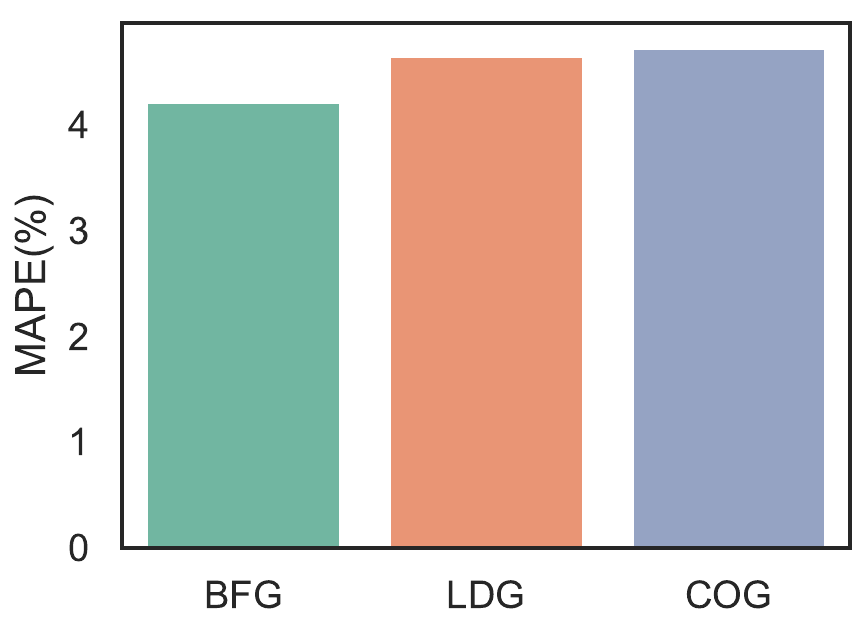}
	\caption{MAPE metric of the test set.}
	\label{fig:mape}
\end{figure*}

\begin{figure*}[htbp!]
	\centering
	\includegraphics[width=0.7\textwidth]{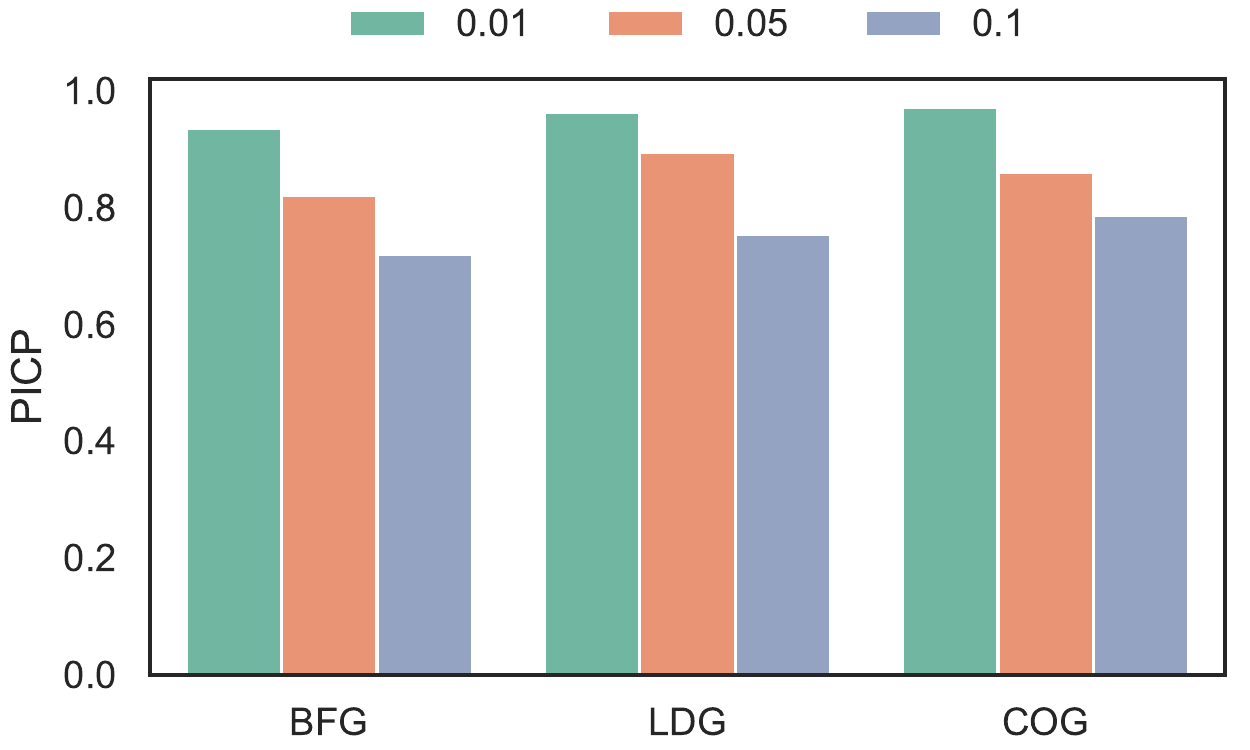}
	\caption{PICP metric of the test set.}
	\label{fig:picp}
\end{figure*}

The calculated metrics are shown in Figure \ref{fig:mape} and \ref{fig:picp}, respectively. Results indicate that: (1) all values of MAPE are less than 0.05, indicating that their nominal predictions are not completely accurate and uncertainty exists. (2) The values of PICP are worse as the quantile declines, hence the quantile selection is important for the proposed DD-TSRO model. Using the validation set (see the attached data file), foretasted nominal values $z_{a,t}^{\circ}$ and deviations $z_{a,t}^{-}$, $z_{a,t}^{+}$ were obtained within eight periods (see the supplemental material Table S7), as shown in Figure \ref{fig:intval}. These values were then used them to parameterize budget-based uncertainty set $\mathcal{Z}$ in the model. 


\begin{figure*}[htbp!]
	\centering
	\includegraphics[width=1.1\textwidth]{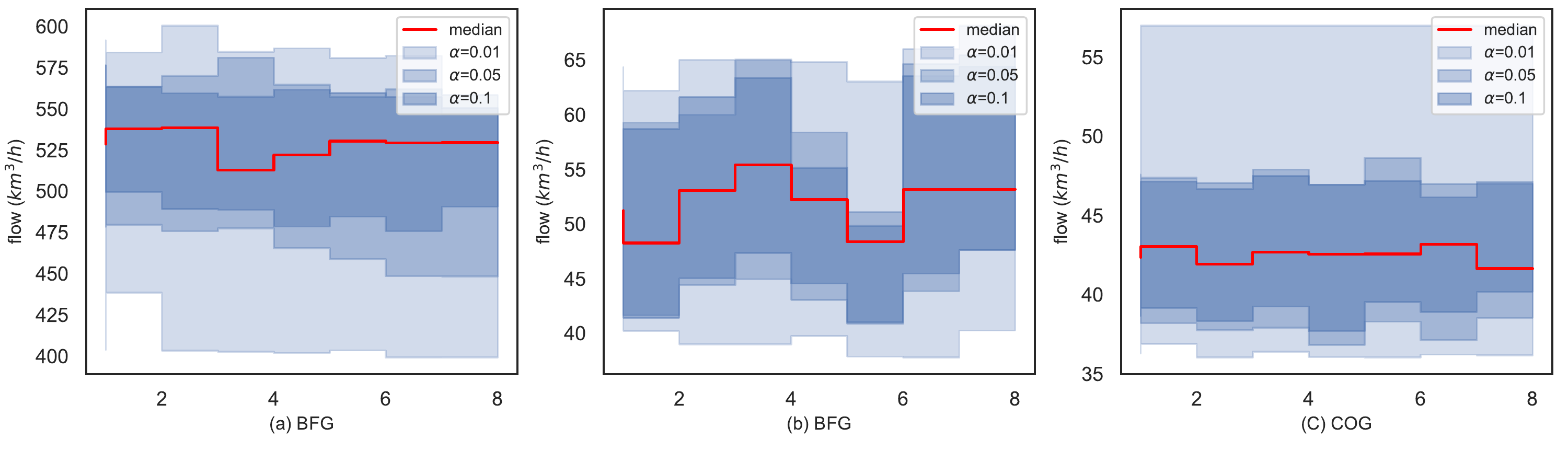}
	\caption{Prediction intervals of the byproduct gases.}
	\label{fig:intval}
\end{figure*}

\subsection{Budget vs. Robustness}

As discussed in Section \ref{S3:model}, the budget parameters $\Gamma_{a}$ were introduced in the TSRO to control the trade-off between the probability of constraint violation and the optimal objective. In this subsection, a set of scenarios were generated to observe the relationship between the budget and the optimal objective under different quantiles. Then the appropriate combination of ($\alpha, \Gamma_a$) was determined to avoid over-conservatism, with the budget $\Gamma_a$ of each gas ranging from 0 to 8. When $\Gamma=0$, it is defined as a deterministic optimization that is the most optimistic without any uncertainty. When $\Gamma=8$, it is defined as a TSRO with a box-based set, which is the most pessimistic. 

\begin{figure*}[htbp!]
	\centering
	\includegraphics[width=0.7\textwidth]{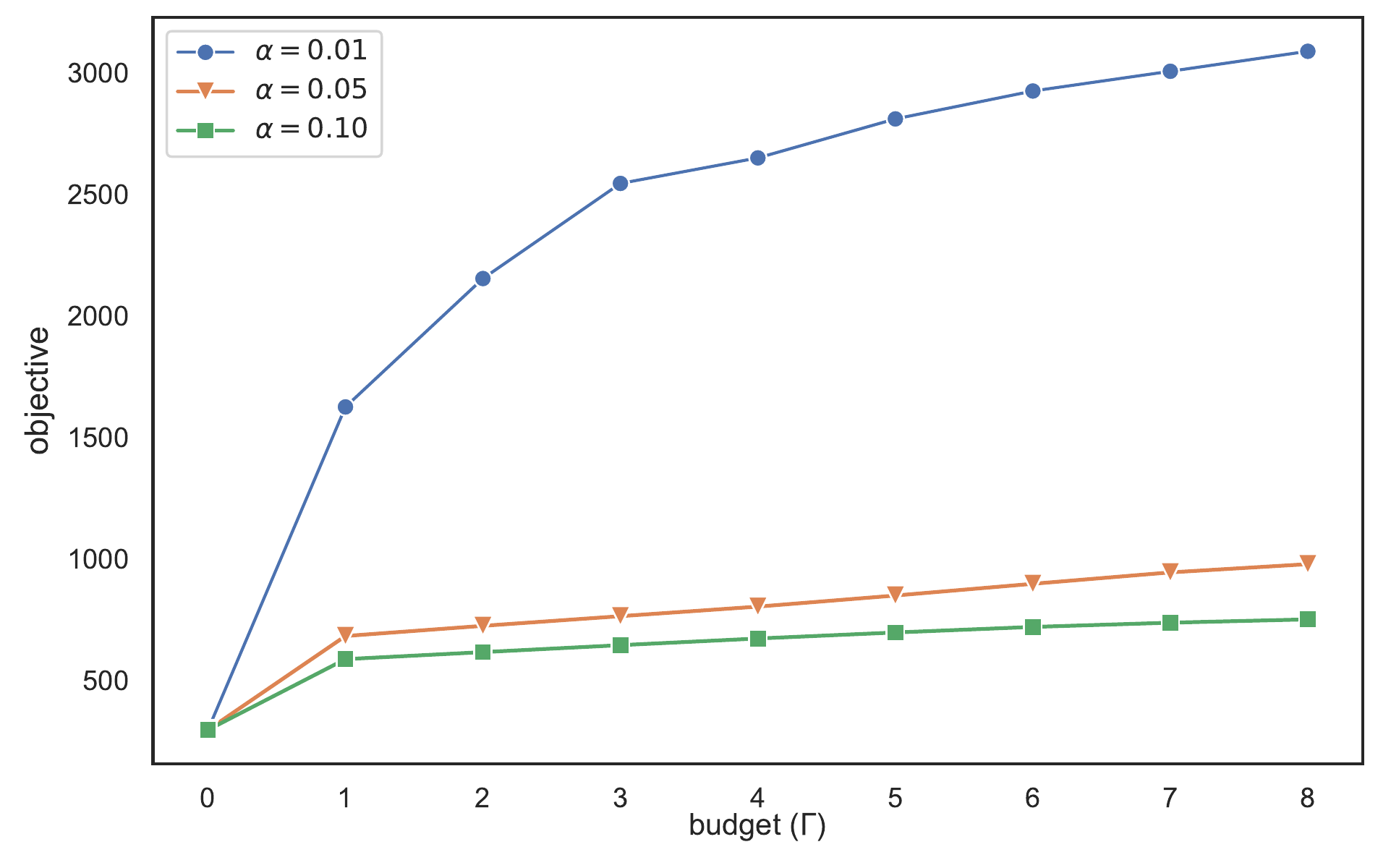}
	\caption{Objective under different quantiles.}
	\label{fig:budget}
\end{figure*}

The variation of optimized objectives is shown in Figure \ref{fig:budget}.  With the increase of the selected budget, the optimal objectives under all quantiles got worse. When $\alpha=0.01$, the optimal objectives were the worst and deteriorated (increase) at the most rapid rate. To trade-off between optimality and robustness, $\alpha=0.05$ and $\Gamma_a =4$ were chosen in the following experiments. 

\subsection{Flexibility vs. Robustness}

In the studied model, the maximum deviation of each gasholder ($\Delta_{k}$)  and the minimal output ratio of each conversion unit ($\eta_{k}^-$) reflect the flexibility of the byproduct gas system as they absorb some uncertainties. The scaling factors of $\Delta_{k}$ were varied from 0.5 to 2.0 with a step of 0.1, whilst $\eta_{k}^-$ was varied from 0.0 to 0.3 with a step of 0.05. 

\begin{figure*}[htbp!]
	\centering
	\includegraphics[width=0.6\textwidth]{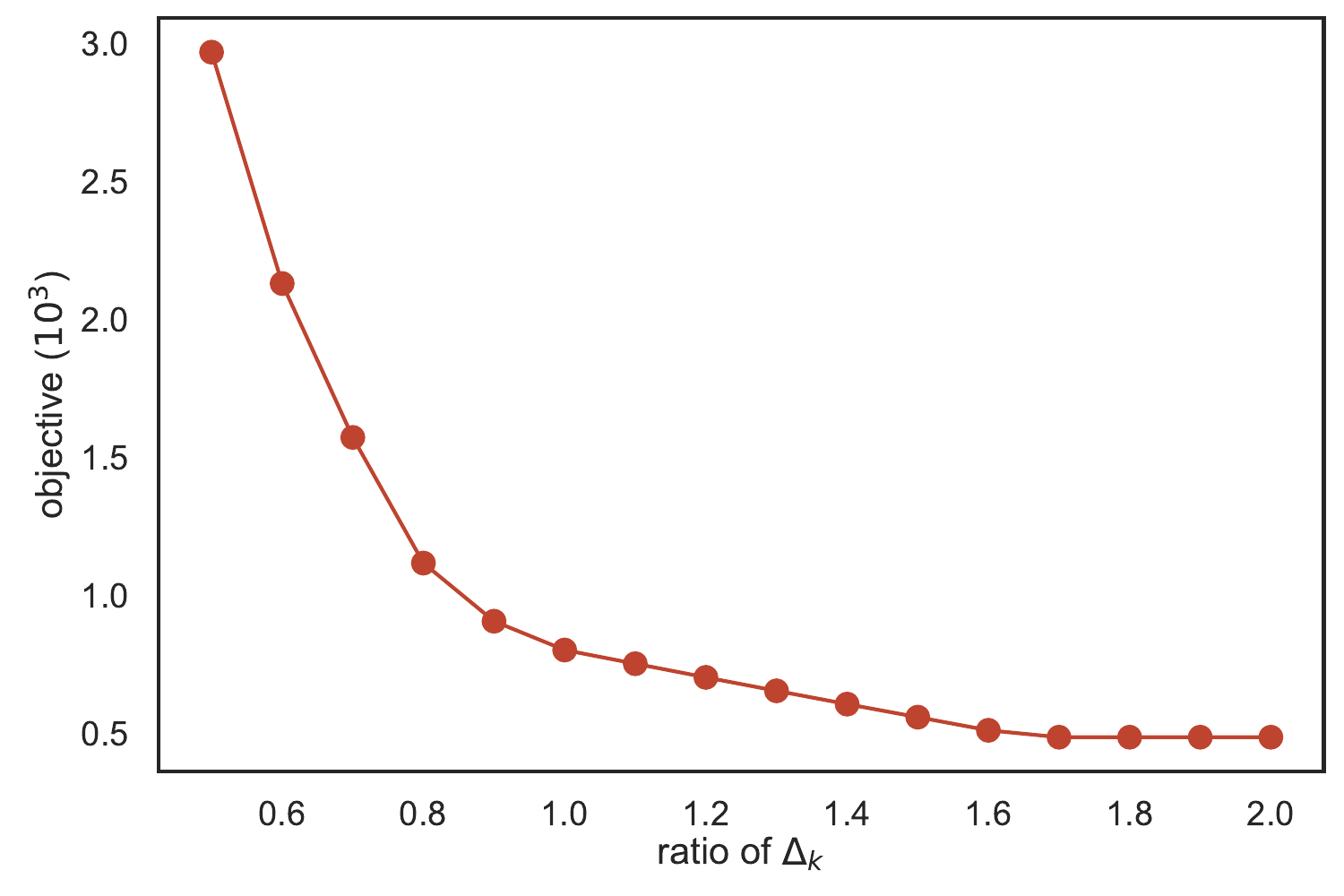}
	\caption{Worst-case objective vs. $\Delta_{k}$.}
	\label{fig:flex1}
\end{figure*}

\begin{figure*}[htbp!]
	\centering
	\includegraphics[width=0.6\textwidth]{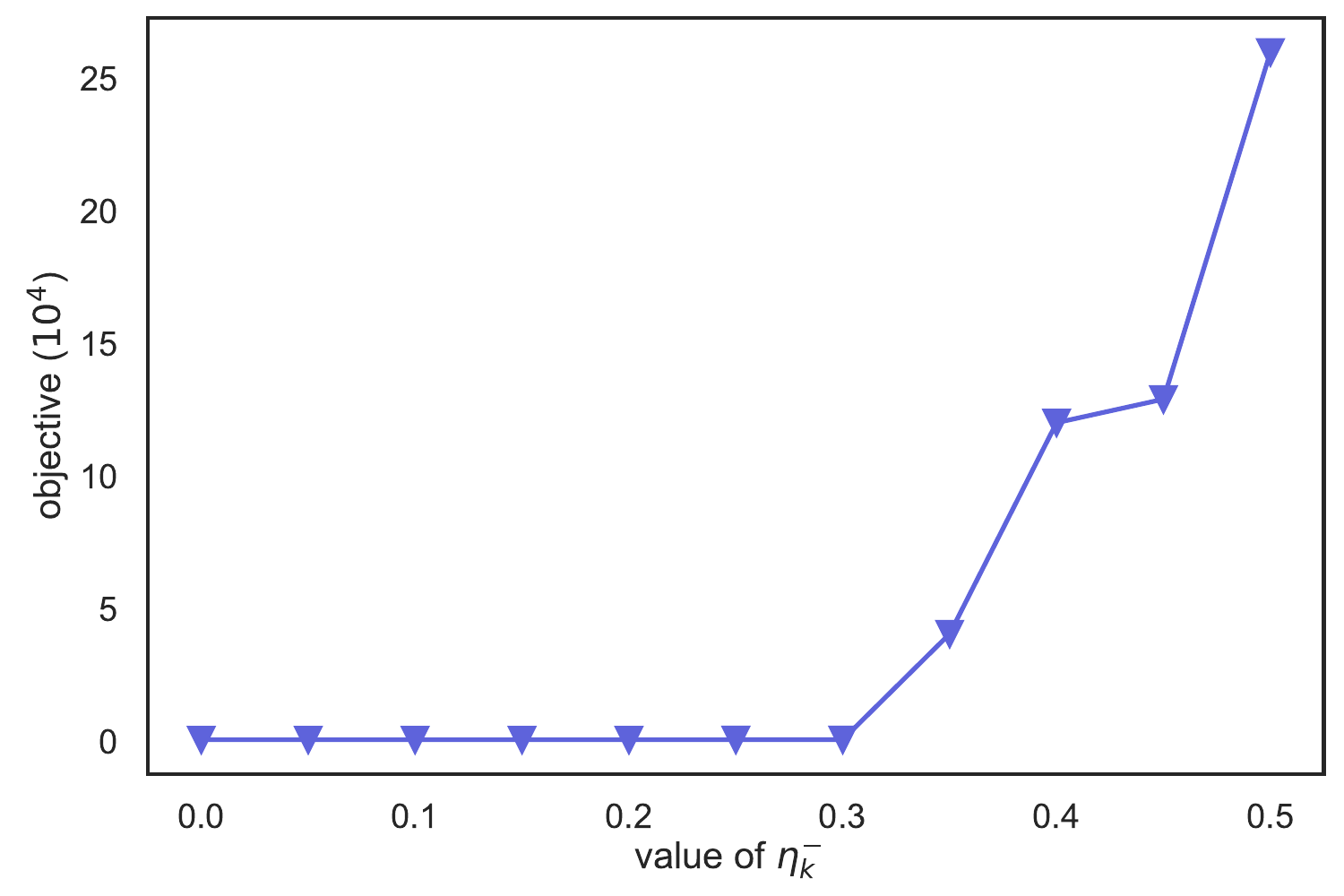}
	\caption{Worst-case objective vs. $\eta_{k}^-$.}
	\label{fig:flex2}
\end{figure*}

The optimized objectives under the worst-case scenario are shown in Figures \ref{fig:flex1} and \ref{fig:flex2}. The curves confirm that (1), with the increase of $\Delta_{k}$, the energy system is more flexible because the optimal objective is increasingly smaller; and (2), when $\eta_{k}^-\geq 0.3 $, the energy system is highly inflexible because optimized objectives become excessively large.

\subsection{Solutions of the proposed approach}

\begin{figure}
	\centering
	\subfigure[Best case]{\includegraphics[width=0.70\textwidth]{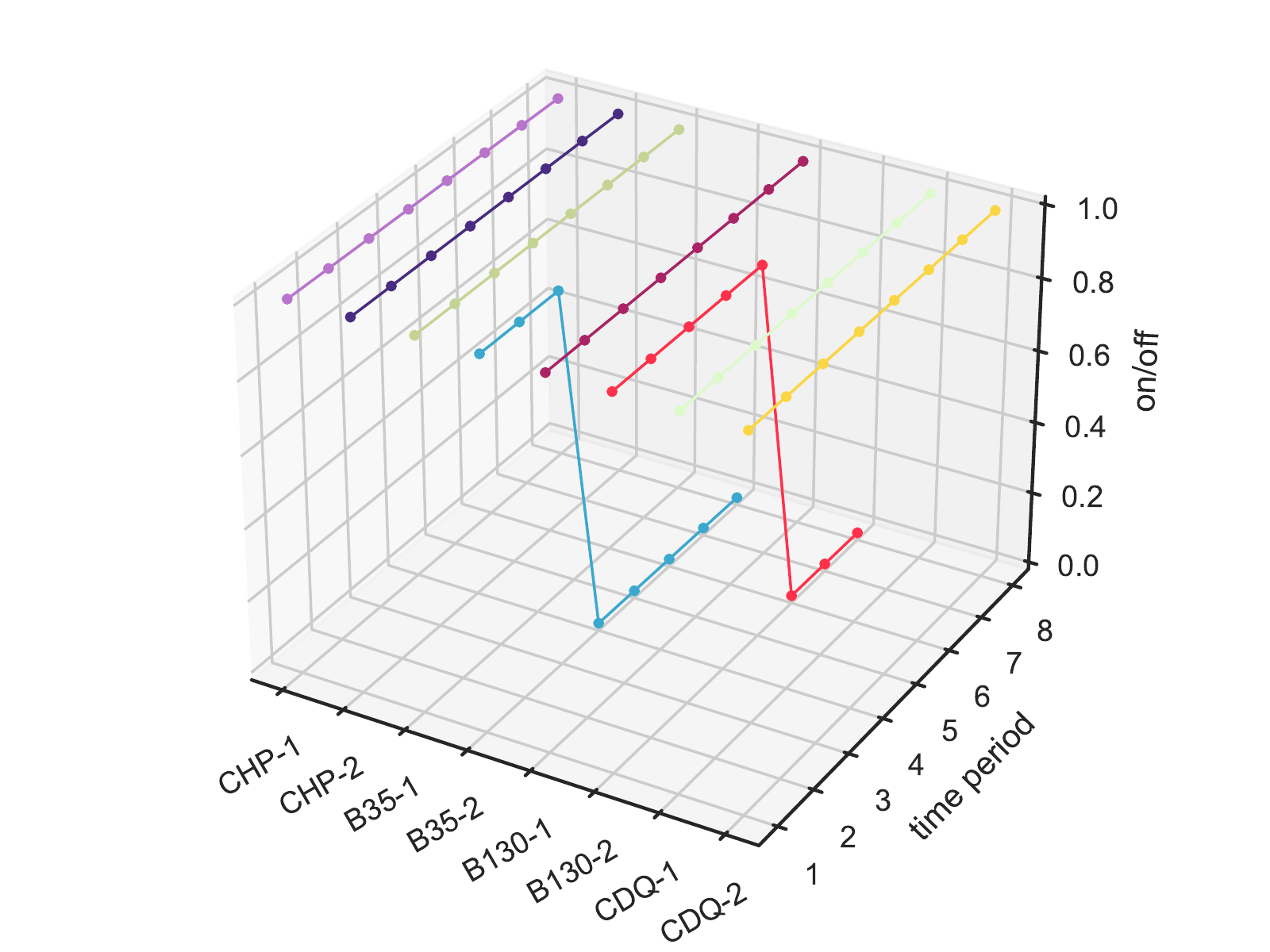}}
	\subfigure[Worst case]{\includegraphics[width=0.70\textwidth]{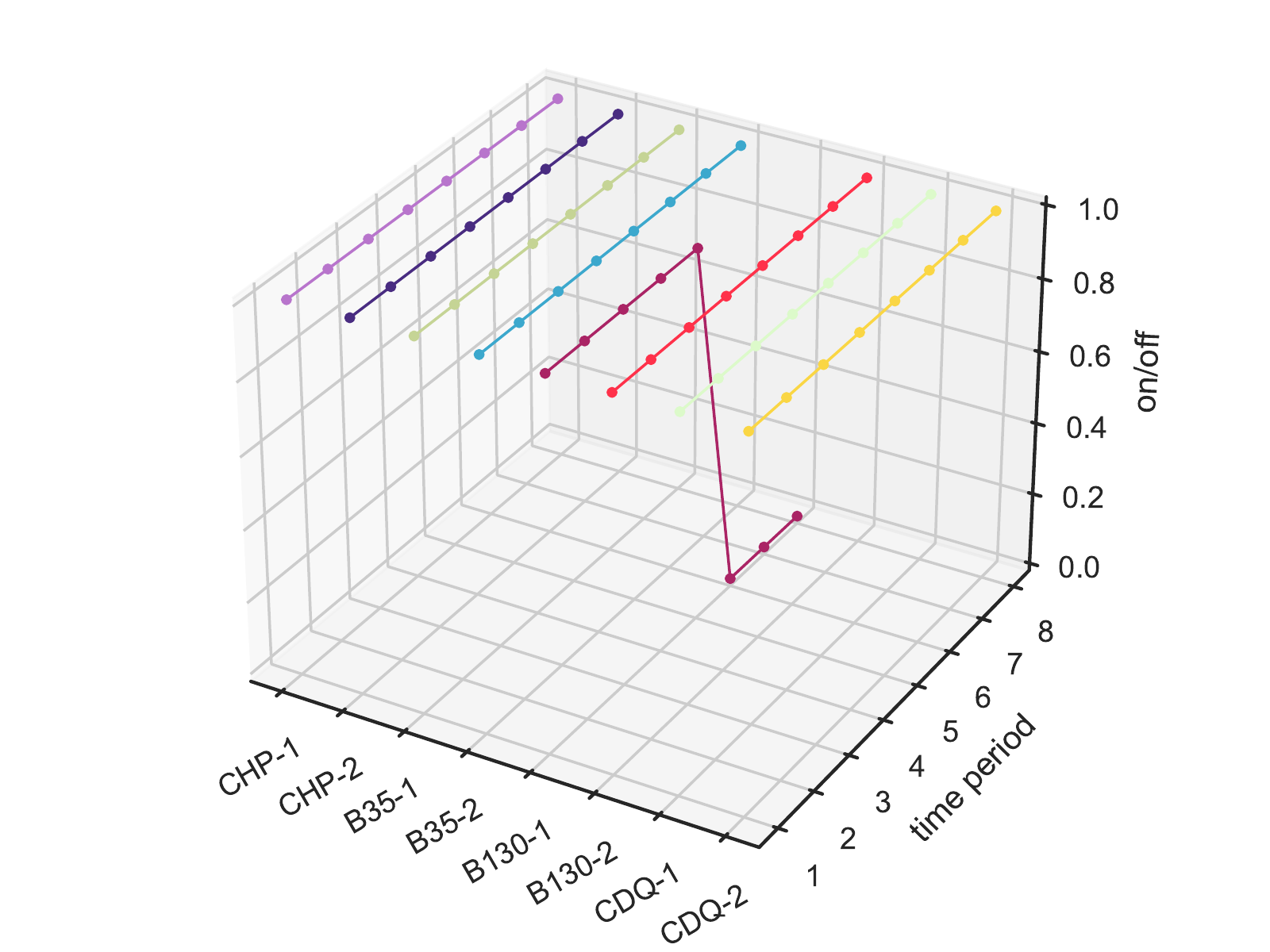}}
	\caption{On-off status of each conversion units.}
	\label{fig:o-var}
\end{figure}

\begin{figure}
	\centering
	\subfigure[Best case]{\includegraphics[width=0.60\textwidth]{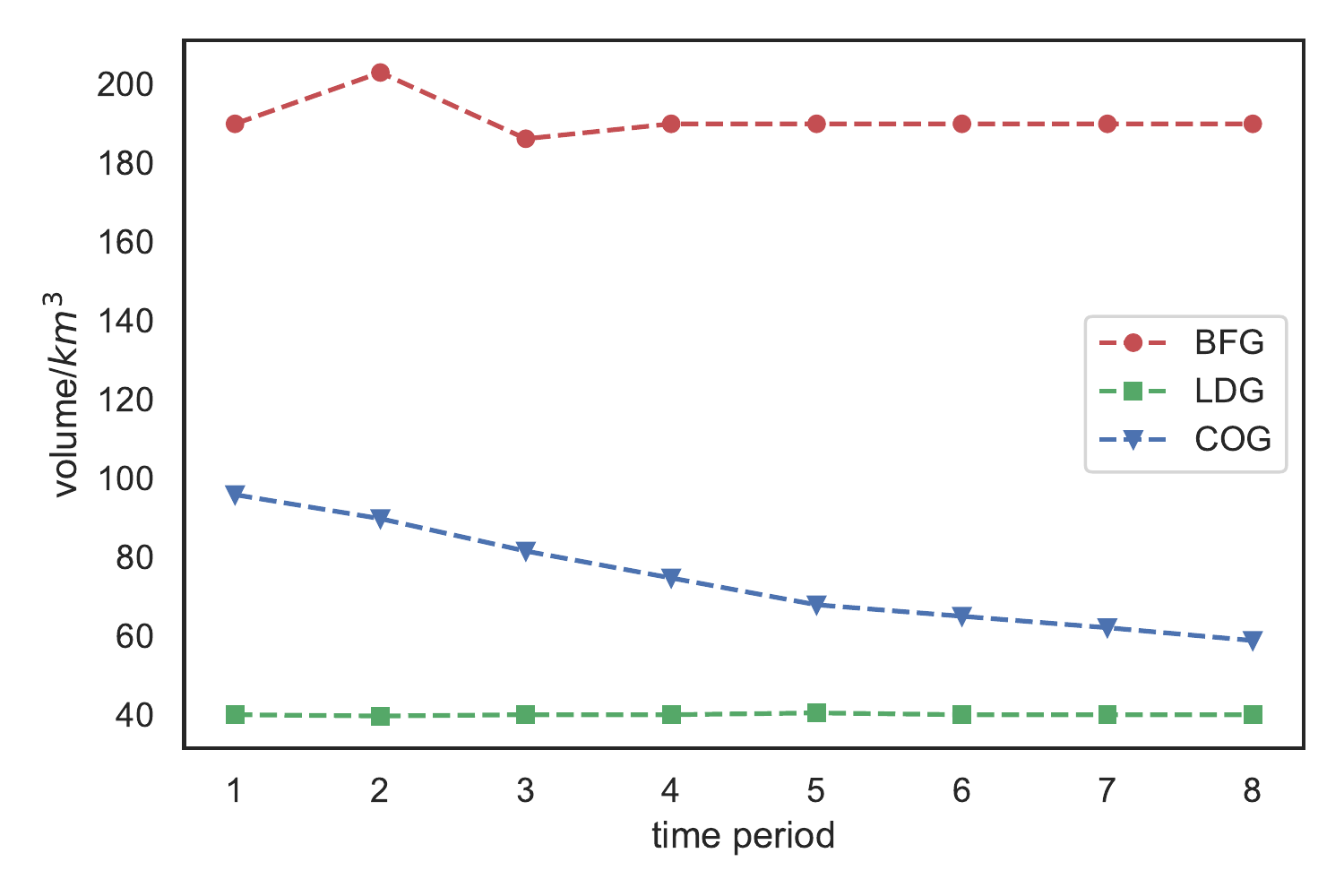}}
	\subfigure[Worst case]{\includegraphics[width=0.60\textwidth]{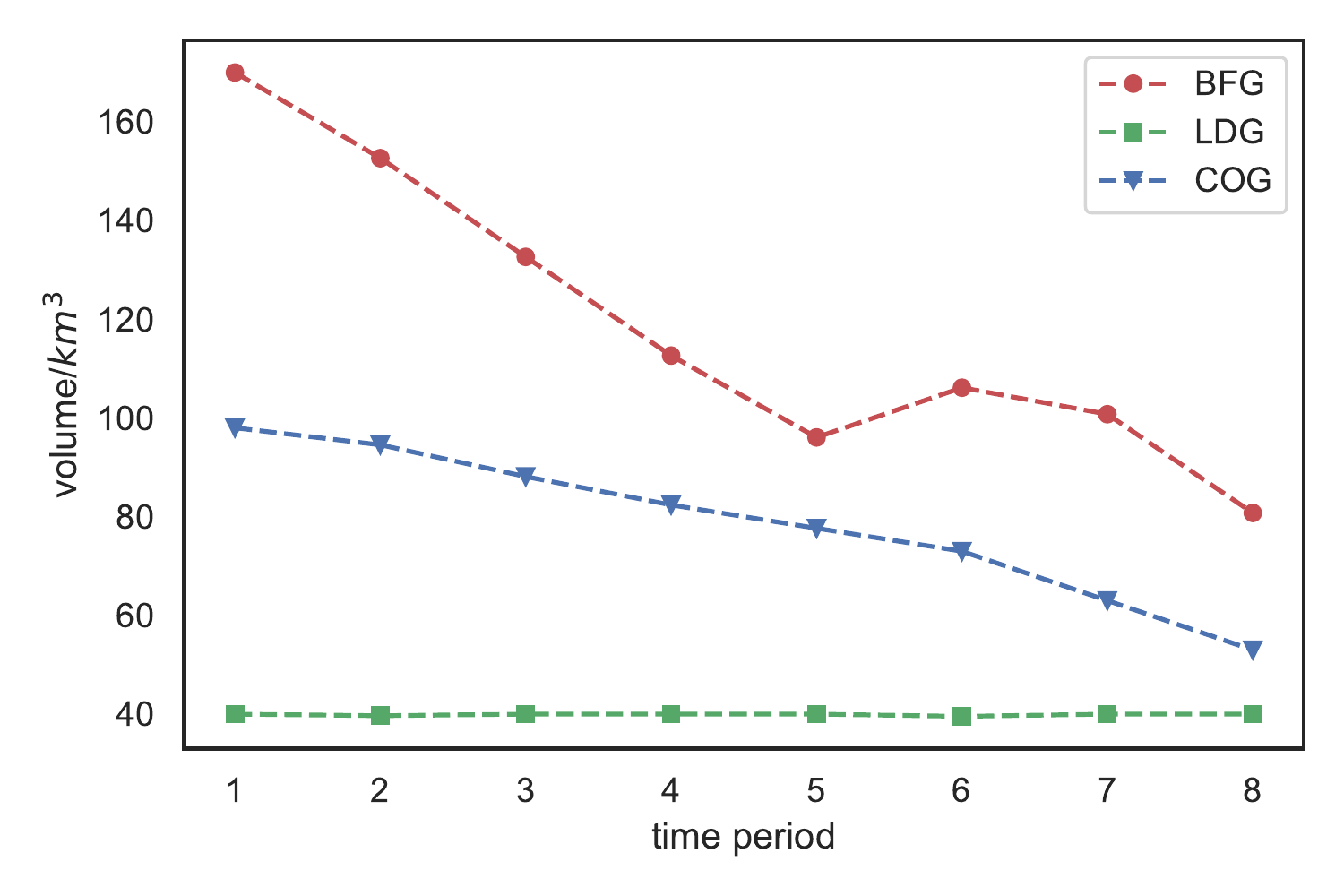}}
	\caption{Volume of each gas holder.}
	\label{fig:u-var}
\end{figure}

\begin{figure}
	\centering
	\subfigure[Best case]{\includegraphics[width=0.70\textwidth]{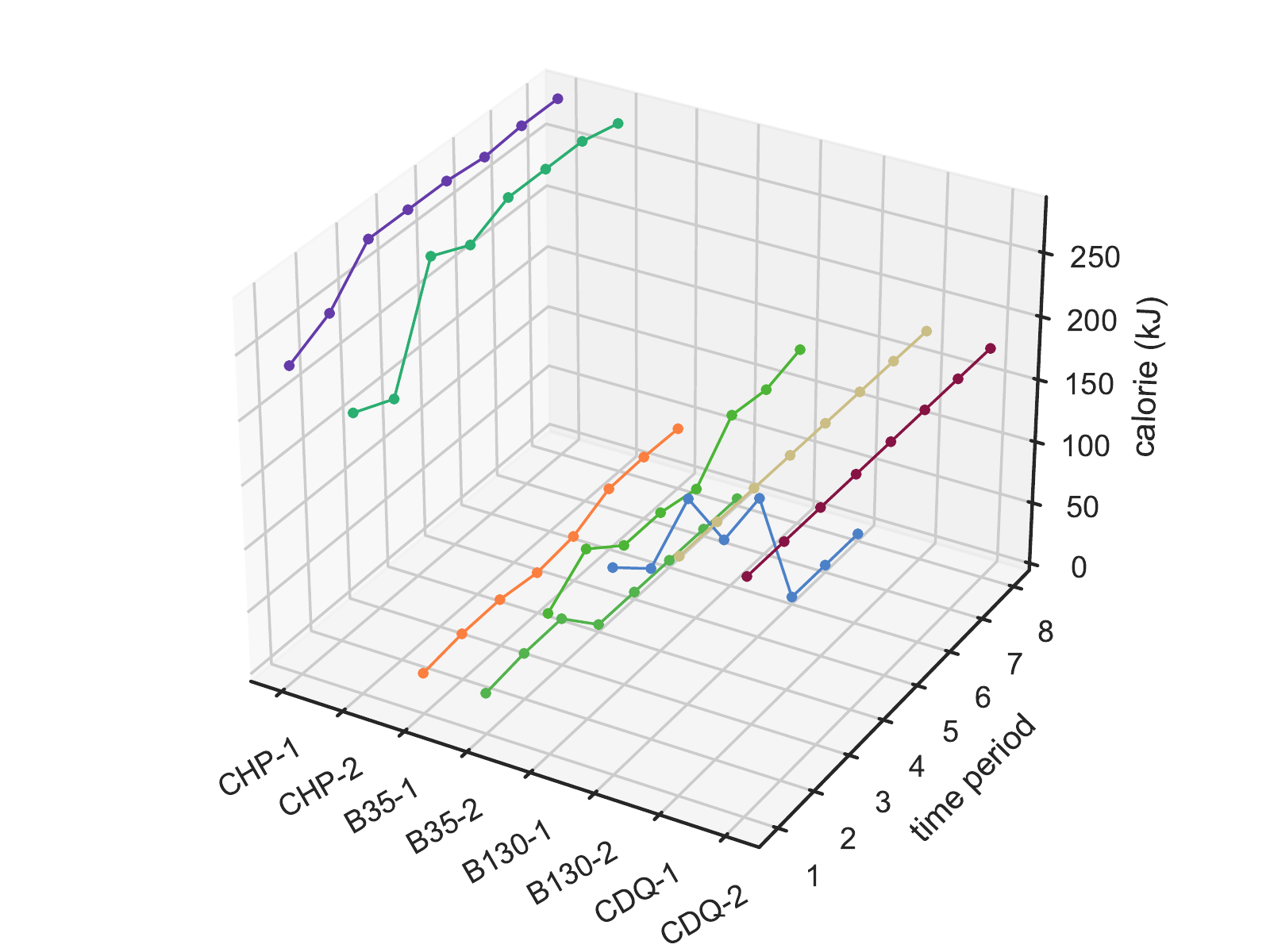}}
	\subfigure[Worst case]{\includegraphics[width=0.70\textwidth]{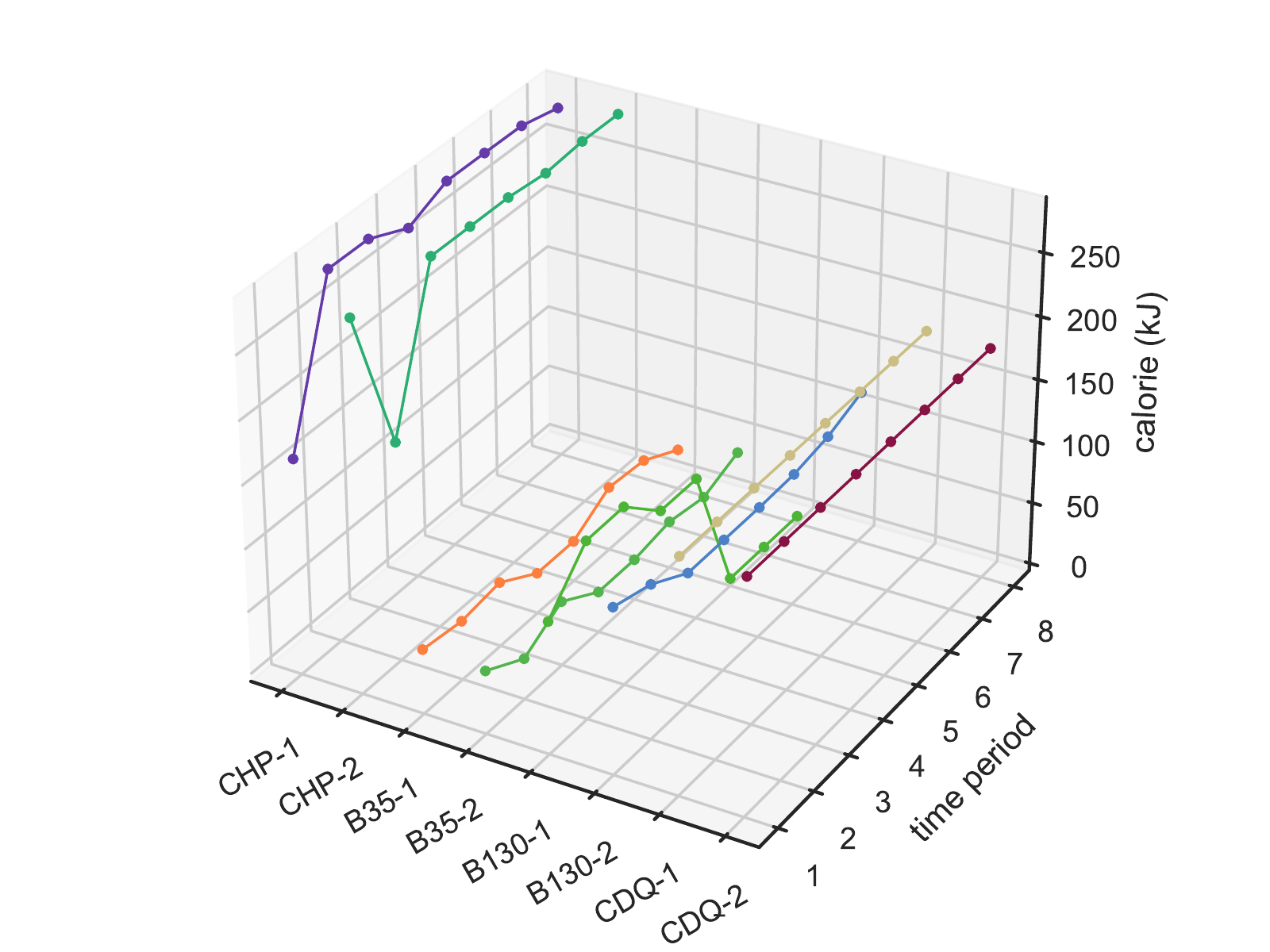}}
	\caption{Input/output calorie of conversion units.}
	\label{fig:f-var}
\end{figure}

With the defined quantile ($\alpha=0.05$), the first-stage decisions under the best-case ($\Gamma_{a}=0$) and the worst-case scenario ($\Gamma_{a}=4$) are shown in Figure \ref{fig:o-var}. The corresponding second-stage decisions under the worst-case scenario are illustrated in Figures \ref{fig:u-var}-\ref{fig:f-var}.

These results suggest that (1), under the worst-case scenario, the capacity of gasholders shrinks since their levels vary more significantly. (2) To avoid potential uncertainty, more conversion units activate because the input and output calories of conversion units also vary significantly. (3) Under both scenarios, all demands are satisfied and no gas is emitted. Thus, the demand outcomes were not plotted.

\section{Conclusions}
\label{S6:concl}

This study investigated the optimal distribution problem of the byproduct gas system under supply uncertainty. Following the supply-storage-conversion-demand network, this study developed an optimal gas distribution model considering uncertain supplies and proposed a TSRO model, including  “here-and-now” decisions, to minimize the start-stop cost of conversion units, and making “wait-and-see” decisions to minimize the operating costs of gasholders and demand penalties. To implement the TSRO model in practice, this work also proposes a quantile regression-based multi-step time series model to quantify the uncertainty of surplus gas, and a column-and-constraint generation algorithm to find the optimal solution. A case study on the industrial energy system of an iron and steel plant was carried out, with results showing that the proposed approach can obtain an effective anvd robust solution:
\begin{enumerate} [(1)]
	\item The quantile-based uncertainty qualification method can forecast both points and intervals.	
	\item The budget under different quantiles makes the trade-off between robustness and optimality.
	\item The flexibility of the storage units and conversion units enables the proposed model to absorb the uncertainty. 
\end{enumerate}

\nomenclature[D01]{$\mathcal{K},k$}{Vertex in the energy network}
\nomenclature[D02]{$\mathcal{A},a$}{Arcs in the energy network}
\nomenclature[D03]{$\mathcal{G},g$}{Energy in the energy network}
\nomenclature[D04]{$\mathcal{T},t$}{Time periods of energy distribution}

\nomenclature[P01]{$\mathcal{A}^{+}(k)$, $\mathcal{A}^{-}(k)$ }{Input and output arcs of vertex $k$}
\nomenclature[P02]{${a}^{+}$, ${a}^{-}$ }{Origin and destination vertex of arc $a$}
\nomenclature[P03]{$e:a$}{Energy of arc $a$}
\nomenclature[P04]{$\omega_{a}$}{Calorie value of the energy flow on arc $a$}
\nomenclature[P05]{$\rho_{k}$}{conversion efficiency of energy unit on vertex $k$}
\nomenclature[P06]{$z_{a,t}$}{Suplus gas at time pertiod $t$}
\nomenclature[P07]{$d_{a,t}$}{Demand of energy eat time pertiod $t$}
\nomenclature[P08]{$\eta_{k}^+$}{Minimum calorie value of the input energies to conversion unit $k$}
\nomenclature[P09]{$\eta_{k}^-$}{Minimum ratio of the ouput energies from conversion unit $k$ to their maximum flows}
\nomenclature[P10]{$\overline{U}_{k}$}{Maximum volume of gasholder $k$}
\nomenclature[P11]{$\underline{U}_{k}$}{Minimum volume of gasholder $k$}
\nomenclature[P12]{$\Delta_{k}$}{Minimum devivation of gasholder $k$}		
\nomenclature[P13]{$\overline{F}^{+}_{k,a}$}{Maximum flow of input energy $a$ of conversion unit $k$}
\nomenclature[P14]{$\underline{F}^{+}_{k,a}$}{Minimum flow of input energy $a$ of conversion unit $k$}	
\nomenclature[P15]{$\overline{F}^{-}_{k,a}$}{Maximum flow of output energy $a$ of conversion unit $k$}
\nomenclature[P16]{$\underline{F}^{-}_{k,a}$}{Maximum flow of output energy $a$ of conversion unit $k$}
\nomenclature[P17]{$\gamma_1$}{Objective coefficient of the start-stop costs}
\nomenclature[P17]{$\gamma_2$}{Objective coefficient of the gasholders operating costs}
\nomenclature[P18]{$\gamma_{3,k}$}{Objective coefficient of surplus and shorages of demand $k$}	

\nomenclature[V1]{$f_{a,t}$}{Continuous variable, which represents the flow of arc $a$ at time period $t$ }
\nomenclature[V2]{$u_{k,t}$}{Continuous variable, which represents stroaage level of unit $v$ at time period $t$ }
\nomenclature[V3]{$v_{g,t}$}{Continuous variable, which represents shortage or over-stock of energy $e$ at time period $t$ }
\nomenclature[V4]{$O_{k,t}$}{Binary variable, which represents the running state (on and off) of conversion unit $k$  at time period $t$}
\nomenclature[V5]{$S_{k,t}$}{Binary variable, which represents the switch state (start up and shut down) of conversion unit $k$  at time period $t$}

\printnomenclature

\section*{Acknowledgment}

This work is supported by the National Natural Science Foundation of China (No. 61873042).




\begin{thebibliography}{00}


\bibitem{Gahm2016} Gahm, Christian, et al. ``Energy-efficient scheduling in manufacturing companies: A review and research framework." European Journal of Operational Research 248.3 (2016): 744-757.
\bibitem{Fan2021} Fan, Zhiyuan, and S. Julio Friedmann. ``Low-carbon production of iron and steel: Technology options, economic assessment, and policy." Joule 5.4 (2021): 829-862.
\bibitem{Ren2021} Ren, Lei, et al. ``A review of CO2 emissions reduction technologies and low-carbon development in the iron and steel industry focusing on China." Renewable and Sustainable Energy Reviews 143 (2021): 110846.
\bibitem{Zhao2017a} Zhao, Xiancong, Hao Bai, and Juxian Hao. ``A review on the optimal scheduling of byproduct gases in steel making industry." Energy Procedia 142 (2017): 2852-2857.
\bibitem{Markland1980} Markland, Robert E. ``Improving fuel utilization in steel mill operations using linear programming." Journal of Operations Management 1.2 (1980): 95-102.
\bibitem{Akimoto1991} Akimoto, K., et al. ``An optimal gas supply for a power plant using a mixed integer programming model." Automatica 27.3 (1991): 513-518.

\bibitem{Kim2003a} Kim, J. H., H-S. Yi, and C. Han. ``A novel MILP model for plantwide multiperiod optimization of byproduct gas supply system in the iron-and steel-making process." Chemical Engineering Research and Design 81.8 (2003): 1015-1025.

\bibitem{Kim2003b} Kim, J. H., et al. ``Plant-wide multiperiod optimal energy resource distribution and byproduct gas holder level control in the iron and steel making process under varying energy demands." Computer Aided Chemical Engineering. Vol. 15. Elsevier, 2003. 882-887.

\bibitem{Kong2010} Kong, Haining, et al. ``An MILP model for optimization of byproduct gases in the integrated iron and steel plant." Applied Energy 87.7 (2010): 2156-2163.

\bibitem{Zhao2015} Zhao, Xiancong, et al. ``A MILP model concerning the optimisation of penalty factors for the short-term distribution of byproduct gases produced in the iron and steel making process." Applied energy 148 (2015): 142-158.

\bibitem{Zhao2017b} Zhao, Xiancong, et al. ``Optimal scheduling of a byproduct gas system in a steel plant considering time-of-use electricity pricing." Applied Energy 195 (2017): 100-113.

\bibitem{Zeng2018} Zeng, Yujiao, et al. ``A novel multi-period mixed-integer linear optimization model for optimal distribution of byproduct gases, steam and power in an iron and steel plant." Energy 143 (2018): 881-899.
\bibitem{Hu2022} Hu, Zhengbiao, and Dongfeng He. ``Operation scheduling optimization of gas–steam–power conversion systems in iron and steel enterprises." Applied Thermal Engineering (2022): 118121.
\bibitem{Zhao2016} Zhao, Jun, et al. ``Data-based predictive optimization for byproduct gas system in steel industry." IEEE Transactions on Automation Science and Engineering 14.4 (2016): 1761-1770.
\bibitem{Jin2018} Jin, Feng, et al. ``A joint scheduling method for multiple byproduct gases in steel industry." Control Engineering Practice 80 (2018): 174-184.
\bibitem{Pena2019} Pena, João G. Coelho, Valter B. de Oliveira Junior, and José L. Félix Salles. ``Optimal scheduling of a by-product gas supply system in the iron-and steel-making process under uncertainties." Computers \& Chemical Engineering 125 (2019): 351-364.

\bibitem{Xi2021} Xi, Han, et al. ``Artificial intelligent based energy scheduling of steel mill gas utilization system towards carbon neutrality." Applied Energy 295 (2021): 117069.

\bibitem{Bertsimas2011} Bertsimas, Dimitris, David B. Brown, and Constantine Caramanis. "Theory and applications of robust optimization." SIAM review 53.3 (2011): 464-501.

\bibitem{Gabrel2014} Gabrel, Virginie, Cécile Murat, and Aurélie Thiele. "Recent advances in robust optimization: An overview." European journal of operational research 235.3 (2014): 471-483.

\bibitem{Rahimian2019} Rahimian, Hamed, and Sanjay Mehrotra. "Distributionally robust optimization: A review." arXiv preprint arXiv:1908.05659 (2019).

\bibitem{Ning2019} Ning, Chao, and Fengqi You. "Optimization under uncertainty in the era of big data and deep learning: When machine learning meets mathematical programming." Computers \& Chemical Engineering 125 (2019): 434-448.

\bibitem{Delage2010} Delage, Erick, and Yinyu Ye. "Distributionally robust optimization under moment uncertainty with application to data-driven problems." Operations research 58.3 (2010): 595-612.

\bibitem{Bertsimas2018} Bertsimas, Dimitris, Vishal Gupta, and Nathan Kallus. "Data-driven robust optimization." Mathematical Programming 167.2 (2018): 235-292.

\bibitem{Shang2017} Shang, Chao, Xiaolin Huang, and Fengqi You. "Data-driven robust optimization based on kernel learning." Computers \& Chemical Engineering 106 (2017): 464-479.

\bibitem{Ning2018} Ning, Chao, and Fengqi You. "Data-driven decision making under uncertainty integrating robust optimization with principal component analysis and kernel smoothing methods." Computers \& Chemical Engineering 112 (2018): 190-210.

\bibitem{Qiu2022} Qiu, Haifeng, et al. "Application of two-stage robust optimization theory in power system scheduling under uncertainties: A review and perspective." Energy (2022): 123942.

\bibitem{Zhao2019}  Zhao, Liang, Chao Ning, and Fengqi You. "A data-driven robust optimization approach to operational optimization of industrial steam systems under uncertainty." Computer Aided Chemical Engineering 46 (2019): 1399-1404.

\bibitem{Shen2020}  Shen, Feifei, et al. "Large-scale industrial energy systems optimization under uncertainty: A data-driven robust optimization approach." Applied Energy 259 (2020): 114199.

\bibitem{Bertsimas2006} Bertsimas, Dimitris, and Aurélie Thiele. "A robust optimization approach to inventory theory." Operations research 54.1 (2006): 150-168.

\bibitem{Mamani2017} Mamani, Hamed, Shima Nassiri, and Michael R. Wagner. "Closed-form solutions for robust inventory management." Management Science 63.5 (2017): 1625-1643.

\bibitem{Koenker1978} Koenker, Roger, and Gilbert Bassett Jr. "Regression quantiles." Econometrica: journal of the Econometric Society (1978): 33-50.

\bibitem{Friedman2001} Friedman, Jerome H. "Greedy function approximation: a gradient boosting machine." Annals of statistics (2001): 1189-1232.

\bibitem{Guo2016} Guo, Yuanxiong, and Chaoyue Zhao. "Islanding-aware robust energy management for microgrids." IEEE Transactions on Smart Grid 9.2 (2016): 1301-1309.

\bibitem{Zeng2013} Zeng, Bo, and Long Zhao. "Solving two-stage robust optimization problems using a column-and-constraint generation method." Operations Research Letters 41.5 (2013): 457-461.

\bibitem{Pedregosa2011} Pedregosa, Fabian, et al. "Scikit-learn: Machine learning in Python." the Journal of machine Learning research 12 (2011): 2825-2830.

\bibitem{Hart2011} Hart, William E., Jean-Paul Watson, and David L. Woodruff. "Pyomo: modeling and solving mathematical programs in Python." Mathematical Programming Computation 3.3 (2011): 219-260.



\end{thebibliography}


\end{document}